\documentclass[12pt]{article}
\usepackage{amssymb}
\usepackage{amsfonts}
\usepackage{graphicx}
\usepackage{amsmath}
\title{\Large Sharpness for $C^1$ linearization of planar
\\
hyperbolic diffeomorphisms
\thanks{Supported by NSFC and MOE research grants}
}
\author{
{\bf Wenmeng Zhang}\,$^{a,b}$~~~and~~~ {\bf Weinian Zhang}\,$^{a}$
\thanks{Corresponding to: matzwn@126.com or matwnzhang@yahoo.com.cn}
\\
$^{a}${\small Yangtze Center of Mathematics and Department of Mathematics}\\
{\small Sichuan University, Chengdu, Sichuan 610064, P. R. China}
\\
$^{b}${\small College of Mathematics Science, Chongqing Normal University}\\
{\small Chongqing 400047, P. R. China}
}

\date{}
\begin{document}
\maketitle

%%%%%%%%%%%%%%%%
\begin{abstract}
Planar hyperbolic diffeomorphisms can be referred to two cases: Poincar\'{e} domain (both eigenvalues lie inside
the unit circle $S^1$) and Siegel domain (one eigenvalue inside $S^1$ but the other outside $S^1$).
In Poincar\'{e} domain it was proved that
$C^{1,\alpha}$ smoothness with $\alpha_0:=1-\log|\lambda_2|/\log|\lambda_1|<\alpha\le 1$, where $\lambda_1$ and $\lambda_2$
are both eigenvalues such that $0<|\lambda_1|<|\lambda_2|<1$,
admits $C^1$ linearization
and the linearization is actually $C^{1,\beta}$.
While a sharp H\"older exponent $\beta>0$ is given,
an interesting problem is: Is the exponent $\alpha_0$ also sharp?
On the other hand, in Siegel domain
we only know that
$C^{1,\alpha}$ smoothness with $\alpha\in (0,1]$ admits $C^1$ linearization.
In this paper we further study the sharpness for $C^1$ linearization in both cases.

\vspace{0.2cm}


\end{abstract}
%%%%%%%%%%%%%%%%

\renewcommand{\theequation}{\thesection.\arabic{equation}}
\newtheorem{lm}{Lemma}
\newtheorem{thm}{Theorem}
\newtheorem{cor}{Corollary}
\newtheorem{fa}{Fact}
\newtheorem{re}{Remark}
\allowdisplaybreaks

\vskip 0.4cm
\parskip 10pt

%%%%%%%%%%%%%%%%%%%%%%%%%%%%%%%%%%%%%%%%%%%%%%%%%%SS1

%%%%%%%%%%%%%%%%%%%%%%%%%%%%%%%%%%%%%%%%%%%%%%%%%%SS1

\section{Introduction}
\setcounter{equation}{0}

Let $(X,\|\cdot\|)$ be a Banach space and
$F:X\to X$ be a diffeomorphism such that
\begin{eqnarray}
F(O)=O
\quad {\rm and} \quad
DF(O)=\Lambda,
\label{def-FR2}
\end{eqnarray}
where $O$ is the origin and $DF(O)$ is the (Fr\'{e}chet)
differentiation of $F$ at $O$. Thus $\Lambda$ is a bounded linear
operator defined on $X$. The local $C^r$ linearization of $F$ is to
find a $C^r$ diffeomorphism $\Phi$ near $O$ such that the conjugacy
equation
\begin{eqnarray}
\Phi\circ F=\Lambda\circ \Phi
\label{schd-eqn}
\end{eqnarray}
holds. The well known Hartman-Grobman Theorem
(\cite{Hart64,pugh-AJM69}) says that $C^1$ diffeomorphisms on $X$
can be $C^0$ linearized near hyperbolic fixed points. Here a fixed
point of $F$ is said to be {\it hyperbolic} if $\Lambda$ has no
eigenvalues on the unit circle $S^1$. In order to preserve more
dynamical properties in the procedure of linearization, one expects
the solution $\Phi$ of equation (\ref{schd-eqn}) to be as regular as
possible. This work goes back to Poincar\'{e} ({\cite{Poincare}}),
who investigated analytic linearization for analytic
diffeomorphisms. Results on $C^r$ linearization for $C^{k}$
diffeomorphisms with $1\le r\le k\le \infty$, initiated by Sternberg
(\cite{ster57,ster58}) in 1950's, can be found in
\cite{Beli-RMS78,Bronkopa-book94,Sell-AJM85}.

$C^1$ linearization is of special interest because it preserves smooth dynamical behaviors
and distinguishes characteristic directions of the systems.
Its applications can be referred to \cite{AAIS-book94,DBo-JDE89} for homoclinic
bifurcations, \cite{FMT-Ann07} for stability of topological mixing of hyperbolic flows,
\cite{HollMelb-07} for Lorenz attractors, \cite{HomHow-CMP99} for Homoclinic tangencies,
and \cite{ZZJZ-JMAA12} for $C^1$ iterative roots of mappings.
For these reasons great efforts (see e.g. \cite{El-JFA01,
R-S-JDE04,R-S-JDDE04} and referneces therein) have been made to
$C^1$ linearization of hyperbolic diffeomorphisms in Euclidean
spaces and Banach spaces since Hartman's \cite{Hart60} and
Belitskii's \cite{Beli-RMS78}. Noting some examples (see
\cite[p.139]{Kuczma1968} and \cite{ster-Duke})
of 1-dimensional $C^1$ hyperbolic mappings which cannot be $C^1$
linearized, one usually considers $C^{1,\alpha}$ mappings with
$\alpha\in (0,1]$, where $C^{1,\alpha}$ denotes the class of all
$C^1$ mappings $F$ whose derivatives satisfy
\begin{eqnarray}
\sup_{x\ne y}\frac{\|D F(x)-D F(y)\|}{\|x-y\|^{\alpha}}<\infty.
\label{Cnorm}
\end{eqnarray}
In spite of some more results on $C^1$ linearization of 1-dimensional mappings (see Theorems 6.2 and 6.3 in \cite{Kuczma1968}),
an important conclusion is that 1-dimensional $C^{1,\alpha}$ hyperbolic mappings can be $C^{1,\alpha}$ linearized for all $\alpha\in (0,1]$,
a corollary of Theorem 6.1 in \cite{Kuczma1968}.
More attentions are paid to 2-dimensional or higher-dimensional $C^{1,\alpha}$ mappings.
Let $F$ be a planar hyperbolic mapping, $X=\mathbb{R}^2$ and
$\lambda:=(\lambda_1,\lambda_2)$, where $\lambda_1$ and $\lambda_2$
are eigenvalues of the linear part $\Lambda$. As indicated in
\cite{Arnold83}, there are two cases in discussion: $\lambda$ lies
in the {\it Poincar\'{e} domain} (i.e., either $\lambda_1$ and
$\lambda_2$ both lie inside the unit circle $S^1$ or both outside
$S^1$); $\lambda$ lies in the {\it Siegel domain} (i.e., the
complement of the Poincar\'{e} domain). In Poincar\'{e} domain, it
suffices to discuss in the case $0<|\lambda_1|\le |\lambda_2|<1$
because the case of expansion can be reduced to this case by
considering the inverse of the mapping.
It is known from \cite[Corollary 1.3.3]{Chap-TMIS02} that
all $C^{1,\alpha}$ mappings can be $C^{1,\alpha}$ linearized if
\begin{eqnarray}
\alpha>\alpha_1:=\log|\lambda_1|/\log|\lambda_2|-1.
\label{a01}
\end{eqnarray}
Since $\alpha_1=0$ when $0<|\lambda_1|=|\lambda_2|<1$, i.e.,
such a $C^{1,\alpha}$ linearization holds for all $\alpha\in (0,1]$
in the particular case, further efforts were made only for the generic case that
\begin{eqnarray}
0<|\lambda_1|<|\lambda_2|<1.
\label{eig-val}
\end{eqnarray}
Recently, the authors (\cite{WZWZ-JFA01}) straightened up the
invariant manifold tangent to the weaker contractive eigen-direction
so as to estimate the convergence rate of the sequence defined by
differentiating the projection of iterate $F^n$ of $F$ onto
the stronger contractive eigen-direction more precisely by $(|\lambda_1|^n)_{n\in\mathbb{N}}$, which guarantees
the convergence of the sequence
$(\Lambda^{-n}F^n)_{n\in\mathbb{N}}$ in $C^1$ norm and gives a $C^1$ solution of equation (\ref{schd-eqn}),
and therefore relaxed the bound $\alpha_1$ to
\begin{eqnarray}
\alpha_0:= 1-\log|\lambda_2|/\log |\lambda_1|
\label{a00}
\end{eqnarray}
for $C^1$ linearization. Moreover, it showed in \cite{WZWZ-JFA01} that
the linearization is actually $C^{1,\beta}$ with a constant $\beta> 0$.
On the other hand, in Siegel domain we assume that
\begin{eqnarray}
0<|\lambda_1|<1<|\lambda_2|.
\label{sdcs}
\end{eqnarray}
It is indicated in the book \cite{AraBeliZhu-book96}
%%by Aranson, Belitskii and Zhuzhoma
that every $C^{1,\alpha}$
diffeomorphism with $\alpha\in(0,1]$ can be $C^1$ linearized.

It is interesting to find sharp bounds for those H\"older exponents $\alpha$ or $\beta$ in these cases.
In 1986 Stowe \cite{Stowe-JDE86} investigated smooth linearization of planar hyperbolic $C^k$ ($k\ge 2$) mappings of the form
$
F(x):=\Lambda x+o(\|x\|^k).
%\label{C2F}
$
In the particular case $k=2$ his results imply $C^{1,\beta}$ linearization
with a sharp estimate of $\beta$.
In the case (\ref{eig-val}), for $\alpha>\alpha_0$ it is proved in \cite[Theorem 2]{WZWZ-JFA01} that
the $C^1$ linearization given in \cite[Theorem 1]{WZWZ-JFA01}
is actaully $C^{1,\beta}$ and estimates of $\beta$ are given as follows:
\begin{description}
\item[{\rm (i)}]
if $\alpha_1<\alpha\le 1$, then $F$ can be $C^{1,\alpha}$ linearized near $O$;

\item[{\rm (ii)}]
if $\alpha=\alpha_1\le 1$, then $F$ can be $C^{1,\beta_1}$ linearized near $O$ for any
$\beta_1\in(0,\alpha)$;

\item[{\rm (iii)}]
if $\alpha=1$ but $\alpha_1> 1$, then $F$ can be $C^{1,\alpha_1^{-1}}$ linearized near $O$;

\item[{\rm (iv)}]
if $\alpha_0<\alpha<\min\{1,\alpha_1\}$, then $F$ can be $C^{1,\beta_2}$ linearized
near $O$, where the exponent
$
\beta_2:=(\alpha_1^{-1}+1)\alpha-1=\alpha_0^{-1}\alpha-1\in (0,1).
$
\end{description}
\vspace{-0.2cm}
A counter example shows that the estimates are sharp.
A natural question is: {\it Is the exponent $\alpha_0$ also sharp?}
In the other case (\ref{sdcs}),
although there is not a question of sharp $\alpha$ and the result of $C^1$ linearization
obtained in \cite{AraBeliZhu-book96} was extended to Banach spaces in
\cite{R-S-JDDE04}, it is still interesting to see if we can
strengthen the $C^1$ linearization to $C^{1,\beta}$ linearization
and give a sharp estimate for $\beta$.

In this paper we further investigate the sharpness for $C^1$ linearization of hyperbolic diffeomorphisms in $\mathbb{R}^2$.
In Poincar\'e domain the remaining case that $\alpha\in (0,\alpha_0]$ will be considered.
We give a counter example in Section 2 to show that
$C^1$ linearization cannot be realized in general,
implying the sharpness of $\alpha_0$,
but rigorously prove in Section 3 a result of
$C^0$ linearization plus differentiability at the fixed point.
In Siegel domain we prove in Section 4 that every $C^{1,\alpha}$ diffeomorphism with $\alpha\in (0,1]$ admits
$C^{1,\beta}$ linearization and give estimates for $\beta$.
It is worthy mentioning that two sequences of mappings corresponding to the projections onto both axes
will be used to approach the desired conjugacy. However, the convergence of the sequences
cannot be proved by using the method of \cite{WZWZ-JFA01}
because divergence factor caused by the part of expansion cannot be avoided.
In order to overcome the difficulty, we not only straighten up both the stable manifold and the unstable one
to make a standard frame, but also
flatten the mapping along these manifolds
for a decomposition. The proof for the decomposition is postponed to Section 5.
Finally, we give counter examples to show that the above-mentioned estimates for $\beta$ are sharp in Section 6.

For convenience, throughout this paper, we assume that
\begin{eqnarray*}
\Lambda:={\rm diag}(\lambda_1,\lambda_2)
\end{eqnarray*}
without loss of generality since $\lambda_1\ne\lambda_2$ in both the
case (\ref{eig-val}) and the case (\ref{sdcs}). For
$x:=(x_1,x_2)\in\mathbb{R}^2$ we consider the norm $\|\cdot\|$
defined by $\|x\|:=\max\{|x_1|,|x_2|\}$. Let $U$ and $V$,
$U\subsetneq V\subset \mathbb{R}^2$, be sufficiently small closed
disks centered at $O$
and let $K,L,M$ and $K_i,L_i,M_i$ ($i\in \mathbb{N}$) be positive
constants.

%%%%%%%%%%%%%%%%%%%%%%%%%%%%%%%%%%%%%

%\section{Preliminaries}
%\setcounter{equation}{0}

%%------------------------------------------
\section{Sharpness of $\alpha_0$ for $C^1$ linearization}
\setcounter{equation}{0}

In this section we consider the case that $\lambda$ belongs to the
Poincar\'{e} domain, i.e., (\ref{eig-val}) holds, and
prove that the exponent $\alpha_0$ (defined in (\ref{a00})),
given in \cite{WZWZ-JFA01}, is a sharp lower bound for $C^1$ linearization.
This fact will be proved by construction of a diffeomorphism $F_*$ depending on $\alpha\in (0,1]$,
which is $C^{1,\alpha}$ but cannot be $C^1$ linearized near its a fixed point $O$ if $\alpha\in (0,\alpha_0]$.

For given $\alpha\in (0,1]$, consider the function $p_{\alpha}:\mathbb{R}\to \mathbb{R}$
defined by
\begin{eqnarray}
p_{\alpha}(s):=\left\{
\begin{array}{ll}
s^{1+\alpha}, &~~ s\ge 0,
\\
\vspace{-0.3cm}
\\
0, &~~ s<0,
\end{array}
\right.
\label{def-p}
\end{eqnarray}
and the function
\begin{eqnarray*}
q(t):=\left\{
\begin{array}{ll}
e^{\frac{1}{t(t-1)}}, &~~ 0<t<1,
\\
\vspace{-0.3cm}
\\
0, &~~ \mbox {other}.
\end{array}
\right.
%\label{def-qq}
\end{eqnarray*}
Define the function
$u:\mathbb{R}^2\backslash\{O\}\to \mathbb{R}$ as
\begin{eqnarray*}
u(x):=\left\{
\begin{array}{ll}
\int_{-\infty}^{x_1/|x_2|}q(t)dt/\int_{-\infty}^{\infty}q(t)dt, &~~
x_2\ne 0,
\\
0, &~~ x_2=0,~x_1<0,
\\
1, &~~ x_2=0,~x_1>0.
\end{array}
\right.
%\label{def-u}
\end{eqnarray*}
Clearly, $p_{\alpha}$
%given in (\ref{def-p})
is $C^{1,\alpha}$ on $\mathbb{R}$ and $u$
%given in (\ref{def-u})
is $C^\infty$ on $\mathbb{R}^2\backslash \{O\}$ such that
\begin{description}
\item[{\rm (U1)}]
${u}(x_1,x_2)=1$ if $x_1\ge|x_2|$ and ${u}(x_1,x_2)=0$ if $x_1\le
0$,

\item[{\rm (U2)}]
$Du(0,x_2)=0$ for all $x_2\in\mathbb{R}\backslash \{0\}$, and

\item[{\rm (U3)}]
%$\partial_{x_1}{u}(x) \ge 0$ and
$\|D^r{u}(x)\|\,\|x\|^{r}\le K$
for all $r=0,1,2$ and all $x\in \mathbb{R}^2\backslash \{O\}$.
\end{description}
Thus we can define a planar mapping $F_*:\mathbb{R}^2\to \mathbb{R}^2$ as
\begin{eqnarray}
F_*(x):=\left\{
\begin{array}{ll}
(\lambda_1x_1+\varrho(x)\,{u}(\lambda_1x_1,\lambda_2x_2)\,
p_\alpha(\lambda_2x_2), \lambda_2x_2), & x\in \mathbb{R}^2\backslash \{O\},
\\
\vspace{-0.3cm}
\\
O, & x=O,
\end{array}
\right. ~~~~~
\label{F-def}
\end{eqnarray}
where $\varrho:\mathbb{R}^2\to \mathbb{R}^2$ is a $C^\infty$ function
such that $\varrho(x)=1$ for
all $x\in U$, $\varrho(x)=0$ for all $x\in \mathbb{R}^2\backslash V$
and $\varrho(x)\ge 0$ for all $x\in\mathbb{R}^2$. Such a function $\varrho$ is usually called a {\it bump function}.
Obviously $F_*$ is $C^1$ in $U\backslash \{O\}$. Moreover, $F_*$ is actually a
$C^{1,\alpha}$ diffeomorphism in $U$ with $DF_*(O)=\Lambda$ due to
the following lemma:

\begin{lm}
Let $F_*$ be given in {\rm (\ref{F-def})}. Then
\begin{eqnarray}
&&\|F_*(x)-\Lambda x\|=o(\|x\|), ~~~~~~~~~ \lim_{x\to
O}DF_*(x)=\Lambda,
\label{1deri}
\\
&&\|DF_*(x)-DF_*(y)\|\le L\|x-y\|^\alpha, \qquad \forall x,y\in U.
\label{hold-cond}
\end{eqnarray}
\end{lm}

{\bf Proof}. Both inequalities given in (\ref{1deri}) can be
verified directly. It implies that $F_*$ is $C^1$ in $U$ with
$DF_*(O)=\Lambda$.

Most of efforts are made to inequality (\ref{hold-cond}).
We only consider the case that
$0<\lambda_1<\lambda_2$ because the other cases can be discussed similarly.
Since (\ref{def-p}) and (U1) imply that the mapping $F_*$
defined in (\ref{F-def}) is linear in the second, third and fourth
quadrants, we only need to prove (\ref{hold-cond}) in the first
quadrant.

First, we claim
that (\ref{hold-cond}) holds in the region
\begin{eqnarray}
U_+:=\{x\in U: 0<\lambda_1x_1<\lambda_2x_2\}.
\label{U+}
\end{eqnarray}
In fact, choose $x=(x_1,x_2)$, $y=(y_1,y_2)$ arbitrarily in $U_+$ and let
$$
\tilde{x}=(\tilde{x}_1,\tilde{x}_2):=(\lambda_1x_1,\lambda_2x_2),~~~
\tilde{y}=(\tilde{y}_1,\tilde{y}_2):=(\lambda_1y_1,\lambda_2y_2).
$$
Clearly, $0<\tilde{x}_1<\tilde{x}_2$ and $0<\tilde{y}_1<\tilde{y}_2$ by (\ref{U+}).
Consider the projections
$$
\pi_1 x:= x_1,~~~~\pi_2 x:=x_2,~~~\forall x=(x_1,x_2)\in \mathbb{R}^2.
$$
It suffices to discuss in the case that
\begin{eqnarray*}
\tilde{x}_2<\tilde{y}_2
%\label{xyorder}
\end{eqnarray*}
because the case $\tilde{x}_2>\tilde{y}_2$ is similar and the case $\tilde{x}_2=\tilde{y}_2$ is simple.
Note that $\varrho(x)= 1$ for all $x\in U$ and that the function $x\mapsto x^\alpha$ is $C^{0,\alpha}$.
Moreover,
\begin{eqnarray*}
\|D^ru(z_1,z_2)\|\,|z_2|^r\le K, \qquad \forall z_1,z_2>0,
\end{eqnarray*}
for $r=1,2$
due to (U1) and (U3).
It follows that
\begin{align}
& |\lambda_1\lambda_2|^{-1}\frac{|(\pi_1 F_*)_{x_1}(x)-(\pi_1
F_*)_{x_1}(y)|}{\|x-y\|^{\alpha}}
     \le\frac{|u_{x_1}(\tilde{x}_1,\tilde{x}_2)\tilde{x}_2^{1+\alpha}
                    -u_{x_1}(\tilde{y}_1,\tilde{y}_2)\tilde{y}_2^{1+\alpha}|}{\|\tilde{x}-\tilde{y}\|^{\alpha}}
\nonumber\\
& \le\frac{|u_{x_1}(\tilde{x}_1,\tilde{x}_2)\tilde{x}_2
                    -u_{x_1}(\tilde{y}_1,\tilde{y}_2)\tilde{y}_2|\,|\tilde{x}_2|^{\alpha}}{\|\tilde{x}-\tilde{y}\|^{\alpha}}
                    +\frac{|u_{x_1}(\tilde{y}_1,\tilde{y}_2)\tilde{y}_2|\,|\tilde{x}_2^{\alpha}-\tilde{y}_2^{\alpha}|}
                           {|\tilde{x}_2-\tilde{y}_2|^{\alpha}}
\nonumber\\
& \le\frac{|u_{x_1}(\tilde{x}_1,\tilde{x}_2)\tilde{x}_2
                    -u_{x_1}(\tilde{y}_1,\tilde{x}_2)\tilde{x}_2|\,|\tilde{x}_2|^{\alpha}}{|\tilde{x}_1-\tilde{y}_1|^{\alpha}}
                    +\frac{|u_{x_1}(\tilde{y}_1,\tilde{x}_2)\tilde{x}_2
                    -u_{x_1}(\tilde{y}_1,\tilde{y}_2)\tilde{y}_2|\,|\tilde{x}_2|^{\alpha}}{|\tilde{x}_2-\tilde{y}_2|^{\alpha}}
                    +M
\nonumber\\
& \le M+|u_{x_1x_1}(\xi_1,\tilde{x}_2)\tilde{x}_2^2|^{\alpha}
                     |u_{x_1}(\tilde{x}_1,\tilde{x}_2)\tilde{x}_2-u_{x_1}(\tilde{y}_1,\tilde{x}_2)\tilde{x}_2|^{1-\alpha}
\nonumber\\
& ~~~+(|u_{x_1x_2}(\tilde{y}_1,\xi_2)\xi_2^2
|^{\alpha}+|u_{x_1}(\tilde{y}_1,\xi_2)\xi_2|^{\alpha})
                    |\tilde{x}_2/\xi_2|^{\alpha} |u_{x_1}(\tilde{y}_1,\tilde{x}_2)\tilde{x}_2
                    -u_{x_1}(\tilde{y}_1,\tilde{y}_2)\tilde{y}_2|^{1-\alpha}
\nonumber\\
& \le K_1:=M+2^{1-\alpha}K+2^{2-\alpha}K,
\label{x1}
\end{align}
where $h_{x_i}:=\partial h/\partial x_i$, $h_{x_ix_j}:=\partial^2h/\partial x_i\partial x_j$,
$\xi_1\in (\tilde{x}_1,\tilde{y}_1)$ or $(\tilde{y}_1,\tilde{x}_1)$ (depending on the comparison between $\tilde{x}_1$ and $\tilde{y}_1$),
and $\xi_2\in (\tilde{x}_2,\tilde{y}_2)$.
%% (since $\tilde{x}_2<\tilde{y}_2$ by (\ref{xyorder})).
Similarly, we can give an analogous estimate for $(\pi_1 F_*)_{x_2}(x)-(\pi_1 F_*)_{x_2}(y)$,
which implies together with the estimate (\ref{x1})
that (\ref{hold-cond}) holds in the region $U_+$,
%%$\{x\in U: 0<\lambda_1x_1<\lambda_2x_2\}$,
i.e., the claimed result is proved.

Since the continuity of $DF_*$ implies the continuity of the
function $(x,y)\mapsto \|DF_*(x)-DF_*(y)\|-L\|x-y\|^\alpha$, the
claimed result implies that (\ref{hold-cond}) also holds in the
closure of $U_+$, i.e., $\{x\in U:
0\le\lambda_1x_1\le\lambda_2x_2\}$. Moreover, (\ref{hold-cond}) is
true in $\{x\in U: 0\le\lambda_2x_2\le\lambda_1x_1\}$ since
$u(\lambda_1x_1,\lambda_2x_2)\equiv 1$ when $0<\lambda_2x_2\le
\lambda_1x_1$ as known in (U1). One can also see that
(\ref{hold-cond}) holds in $\{x\in U: \mbox{either $x_1$ or $x_2$ is
$\le 0$}\}$ because $F_*$ is linear in the second, third and fourth
quadrants. On the other hand, we have
\begin{eqnarray*}
U
&=&\{x\in U:
0\le\lambda_1x_1\le\lambda_2x_2\}\cup \{x\in U:
0\le\lambda_2x_2\le\lambda_1x_1\}
\\
&&\cup ~\{x\in U: \mbox{either $x_1$ or $x_2$ is $\le 0$}\}
\end{eqnarray*}
and
\begin{align*}
\|F_*(x)-F_*(y)\|&\le \|F_*(x)-F_*(z)\|+\|F_*(z)-F_*(y)\|
\\
&\le L_1(\|x-z\|^\alpha+\|z-y\|^\alpha)\le 2L_1\|x-y\|^\alpha,
\end{align*}
where $z$ lies on the line between $x$ and $y$.
Then it follows that inequality
(\ref{hold-cond}) is true in $U$.
The proof is completed. \qquad$\Box$

Before presenting the main result of this section, we need another
useful lemma.

\begin{lm}
Let $F:\mathbb{R}^2\to \mathbb{R}^2$ be $C^{1,\alpha}$
such that {\rm (\ref{def-FR2})} holds and
\begin{eqnarray}
\pi_1F(0,x_2)=0, ~~~~~ \pi_2 F(x)=\lambda_2x_2, ~~~~~\forall x\in U.
\label{F1F2}
\end{eqnarray}
Suppose that {\rm (\ref{eig-val})} holds and that $F$ admits $C^1$
linearization in $U$. Then the limit
\begin{eqnarray}
\widetilde{\Phi}(x):=\lim_{n\to \infty}\Lambda^{-n}F^n(x)
\label{def-Phi22}
\end{eqnarray}
exists uniformly in $U$ and $\widetilde{\Phi}$ is $C^1$.
\label{lm-phiF}
\end{lm}

{\bf Proof}. Since $F$ admits $C^1$ linearization in $U$,
there is a $C^1$ diffeomorphism $\Phi:U\to \mathbb{R}^2$ such that
equation (\ref{schd-eqn}) holds
and therefore
$\Phi(F^n(x))=\Lambda^n\Phi(x)$ for every integer $n\ge 1$.
Differentiating both sides of the equality, we get
\begin{eqnarray}
{D\Phi}(F^n(x))
=\Lambda^n{D\Phi}(x)(DF^n(x))^{-1},~~~ \forall n\in\mathbb{N},~\forall x\in U,
\label{phi'}
\end{eqnarray}
where
$(DF^n(x))^{-1}$ presents the inverse of the matrix $DF^n(x)$.
Knowing from (\ref{phi'}) in the case of $n=1$ that the invertible matrix ${D\Phi}(O)$ commutes
with $\Lambda$, we see from (\ref{eig-val}) that
\begin{eqnarray}
{D\Phi}(O)={\rm diag}(p_1,p_2),
\label{phi0}
\end{eqnarray}
where $p_1,p_2\ne 0$ are both constants.
In virtue of the second equality of (\ref{F1F2}), we
can put
\begin{eqnarray}
DF^n(x):= \left(\begin{array}{cc} {a}_n(x) & {b}_n(x)
\\
0 & \lambda_2^n
\end{array}\right)
~~~ {\rm and} ~~~ {D\Phi}(x):= \left(\begin{array}{cc}
\varphi_{11}(x) & \varphi_{12}(x)
\\
\varphi_{21}(x) & \varphi_{22}(x)
\end{array}\right)
\label{F-phi}
\end{eqnarray}
for each $n\in\mathbb{N}$, where ${a}_n$, ${b}_n$ and
$\varphi_{ij}$'s ($i,j=1,2$) are $C^{0,\alpha}$ functions defined on
$U$. Then, by (\ref{phi'}) and (\ref{F-phi}), we obtain
\begin{eqnarray*}
\varphi_{11}(F^n(x))=\lambda_1^n\varphi_{11}(x)/{a}_n(x), ~~~~~
\forall n\in\mathbb{N}.
%\label{compt-diff}
\end{eqnarray*}
This implies that the limit $\lim_{n\to\infty}\lambda_1^{-n}a_n(x)$
exists uniformly because the sequence $(F^n(x))_{n\in\mathbb{N}}$
converges uniformly to $O$ and $\varphi_{11}(O)=p_1\ne 0$ by
(\ref{phi0}). Therefore
$$
\lim_{n\to\infty}\lambda_1^{-n}{a}_n(x)=\frac{1}{p_1}\varphi_{11}(x),
~~~~~ \forall x\in U.
$$
Note that the first equality of (\ref{F1F2}) implies that $\pi_1
F^n(0,x_2)=0 $ for all $n\in\mathbb{N}$ by induction. Then
\begin{eqnarray}
\lim_{n\to\infty}\lambda_1^{-n}\pi_1F^n(x)
=\lim_{n\to\infty}\int_{0}^{x_1}\lambda_1^{-n}{a}_n(t,x_2)dt
=\frac{1}{p_1}\int_{0}^{x_1}\varphi_{11}(t,x_2)dt
\label{ddd}
\end{eqnarray}
and the limit exists uniformly in $U$.
By (\ref{ddd}), we know that $\widetilde{\Phi}$ given in
(\ref{def-Phi22}) is well defined and
$$
\widetilde{\Phi}(x)= \Bigg(
\frac{1}{p_1}\int_{0}^{x_1}\varphi_{11}(t,x_2)dt,~ x_2 \Bigg), ~~~~~
\forall x\in U.
$$
On the other hand, by the
second equality of (\ref{F-phi}),
$$
\pi_1 \Phi(x)=\int_{0}^{x_1}\varphi_{11}(t,x_2)dt+\varpi(x_2),
$$
where $\varpi$ is a function defined on $U\cap \mathbb{R}$.
Since the statement at the beginning of the proof guarantees that $\pi_1
\Phi(x)$ is $C^1$, putting $x_1=0$, we see that $\varpi$ is $C^1$
and so does the function $\int_{0}^{x_1}\varphi_{11}(t,x_2)dt$.
Finally, we conclude that $\widetilde{\Phi}$ is $C^1$. This completes the proof. \quad$\Box$

Recall that the two constants $\alpha_0$ and $\alpha_1$ are given in
(\ref{a01}) and (\ref{a00}) respectively. Now we are ready to give
the following result:
\begin{thm}
In the case of {\rm (\ref{eig-val})}, the lower bound $\alpha_0$ of $\alpha$ for $C^1$ linearization is
a sharp one.
\label{fa-con}
\end{thm}

{\bf Proof}.
The idea of the proof is to show that the $C^{1,\alpha}$ diffeomorphism $F_*$ defined in (\ref{F-def})
cannot be $C^1$ linearized near its a fixed point $O$ if $\alpha\in (0,\alpha_0]$.
We only need to consider the case that $0<\lambda_1<\lambda_2<1$. Otherwise,
if one of eigenvalues of $DF_*(O)$ is negative, we consider the
quadratic iterate $F_*^2$ instead of $F_*$ to obtain the same
conclusion. Fix a constant $\xi\in U\cap (0,+\infty)$ arbitrarily
and choose
\begin{eqnarray*}
n_0(x_2):=\alpha_1^{-1}\log_{\lambda_2}(x_2/\xi)
\end{eqnarray*}
for all sufficiently small $x_2\in U\cap (0,+\infty)$.
Since the mapping $F_*$
%%defined in (\ref{F-def})
is of the same form
as the one given in \cite[(4.48)]{WZWZ-JFA01} for all $x\in
U$, formula (4.54) of \cite{WZWZ-JFA01} yields
\begin{eqnarray}
\pi_1 F_*^{n}(\xi,x_2)
\ge\lambda_1^{n}(\xi+Mx_2^{(\alpha_1^{-1}+1)\alpha})=
\lambda_1^{n}(\xi+Mx_2^{\alpha_0^{-1}\alpha}), ~~~ \forall n>
n_0(x_2),
\label{pi1F<}
\end{eqnarray}
where $M>0$ is constant and independent of $x_2$ and $n$. Moreover, we
can prove by induction that
\begin{eqnarray}
\pi_1 F_*^{n}(\xi,0)=\lambda_1^n\xi, ~~~~~ \forall n\in\mathbb{N}.
\label{Fxi}
\end{eqnarray}

For a reduction to absurdity, assume that $F_*$
%%given in (\ref{F-def})
admits $C^1$ linearization
in $U$. Obviously, $F_*$ satisfies both equalities given in
(\ref{F1F2}). By Lemma \ref{lm-phiF}, the mapping
\begin{eqnarray}
\Phi_*(x):=\lim_{n\to \infty}\Lambda^{-n}F_*^n(x)
\label{def-Phi*}
\end{eqnarray}
is well defined in $U$ and is $C^1$. Combining (\ref{pi1F<})
with (\ref{Fxi}), we get
\begin{eqnarray}
&&\lim_{x_2\to
0^+}\frac{\pi_1\Phi_*(\xi,x_2)-\pi_1\Phi_*(\xi,0)}{x_2}
\nonumber\\
&&= \lim_{x_2\to
0^+}\frac{\lim_{n\to\infty}\lambda_1^{-n}\pi_1F_*^{n}(\xi,x_2)
    -\lim_{n\to\infty}\lambda_1^{-n}\pi_1 F_*^{n}(\xi,0)}{x_2}
\nonumber\\
&&\ge M\lim_{x_2\to 0^+}x_2^{\alpha_0^{-1}\alpha-1}\ge M>0
\label{derv}
\end{eqnarray}
because $\alpha_0^{-1}\alpha-1\le 0$ provided $\alpha\le\alpha_0$.
This implies that the right derivative of $\pi_1\Phi_*(\xi,x_2)$
with respect to the variable $x_2$ at $0$ is positive.
On the other hand, as mentioned at the beginning of the proof,
$F_*$ is linear, i.e., $F_*(x)=(\lambda_1x_1,\lambda_2x_2)$, in the forth quadrant.
By (\ref{def-Phi*}),
$$
\pi_1\Phi_*(\xi,x_2)=\xi, ~~~~~ \forall x_2\in U\cap(-\infty, 0),
$$
and consequently
\begin{eqnarray}
\lim_{x_2\to
0^-}\frac{\pi_1\Phi_*(\xi,x_2)-\pi_1\Phi_*(\xi,0)}{x_2}=0.
\label{dervleft}
\end{eqnarray}
This implies that the left derivative of $\pi_1\Phi_*(\xi,x_2)$ with
respect to $x_2$ at $0$ is $0$. Hence, the function $\pi_1\Phi_*$ is
not differentiable at the point $(\xi,0)$ by (\ref{derv}) and
(\ref{dervleft}), which contradicts to the fact that $\Phi_*$
is $C^1$, indicated just below (\ref{def-Phi*}). The proof is completed.
\quad$\Box$

%%%%%%%%%%%%%%%%%%%%%%%%%%%%%%%%%%%%%%%%%%%%%%%%%%

\section{Differentiable linearization in Poincar\'{e} domain}
\setcounter{equation}{0}

In this section we continue our discussion in the case of (\ref{eig-val}).
Theorem \ref{fa-con} shows that for $\alpha\in(0,\alpha_0]$ one cannot expect
$C^1$ linearization in general. In this section we persue a weaker result:
$C^0$ linearization plus differentiability at the fixed point.

\begin{lm}
Let $F:\mathbb{R}^2\to \mathbb{R}^2$ be $C^{1,\alpha}$ with $\alpha\in (0,\alpha_0]$ such that
{\rm (\ref{def-FR2})} and {\rm (\ref{eig-val})} hold. Then there is
a $C^{1,\alpha}$ diffeomorphism $\Theta:U\to \mathbb{R}^2$, where
$\Theta(O)=O$ and $D\Theta(O)=$ {\rm id}, the identity mapping in $\mathbb{R}^2$, such that the mapping
$\widetilde{F}:=\Theta\circ F\circ \Theta^{-1}$ is $C^{1,\alpha}$
and satisfies
\begin{align}
\pi_1\widetilde{F}(0,x_2)=0, ~~~
\pi_2\widetilde{F}(x)=\lambda_2x_2,~~~
\forall x_2\in U\cap \mathbb{R},~ \forall x\in U.
\label{flatten}
\end{align}
\label{lm-reduce}
\end{lm}

\vspace{-0.6cm}
{\bf Proof}.
Note that our bump function $\varrho$ satisfies that
$\varrho(x)=1$ for all $x\in U$ and $=0$ for all $x\in
\mathbb{R}^2\backslash V$.
Multiplying the nonlinear part of $F$ by the bump function $\varrho$, we obtain a modified
diffeomorphism, which is still denoted by the same $F$,
such that
\begin{eqnarray}
%F(x)=\Lambda x, \quad \forall x\in \mathbb{R}^2\backslash V,
%\quad {\rm and} \quad
\|DF(x)-\Lambda\|\le \eta, \quad \forall x\in \mathbb{R}^2,
\label{FFbumpR2}
\end{eqnarray}
where $\eta>0$ is a small constant depending on $V$. This
modification does not affect our results at all because we are only
interested in local properties of $F$.

We first claim that the sequence $\big(\lambda_2^{-n}D(\pi_2
F^n)(x)\big)_{n\in\mathbb{N}}$ is uniformly convergent in $U$,
which is proved by the equality
\begin{align}
&\lim_{n\to \infty}\lambda_2^{-n}D(\pi_2 F^n)(x)
\nonumber\\
&=\sum_{i=1}^{\infty}
\Big(\lambda_2^{-(n+1)}D(\pi_2 F^{n+1})(x)-\lambda_2^{-n}D(\pi_2 F^{n})(x)\Big)
+\lambda_2^{-1}D(\pi_2 F)(x).
\label{dfdf}
\end{align}
Note that
\begin{align*}
&|\lambda_2^{-(n+1)}D(\pi_2 F^{n+1})(x)-\lambda_2^{-n}D(\pi_2 F^{n})(x)|
\nonumber\\
&\le |\lambda_2^{-(n+1)}|\,\|D(\pi_2 F^{n+1})(x)-\lambda_2 D(\pi_2 F^{n})(x)\|
\nonumber\\
&\le
|\lambda_2^{-(n+1)}|\,\|DF(F^n(x))-\Lambda\|\,\|DF^n(x)\|
\nonumber\\
&\le K\bigg(\frac{(|\lambda_2|+\eta)^{1+\alpha}}{|\lambda_2|}\bigg)^n,
~~~~~ \forall n\in\mathbb{N},
%\label{bxzl1}
\end{align*}
because $\|DF^n(x)\|\le \|DF(F^{n-1}(x))\|\cdots\|DF(x)\|\le (|\lambda_2|+\eta)^n$, as known from (\ref{FFbumpR2}).
Moreover, $(|\lambda_2|+\eta)^{1+\alpha}/|\lambda_2|<1$ for small $\eta>0$.
It implies the uniform convergence of the series $\sum_{i=1}^\infty (\cdots)$ in (\ref{dfdf}).
Furthermore, one can see that the sequence
$\big(\lambda_2^{-n}\pi_2 F^n(x)\big)_{n\in\mathbb{N}}$ is also
convergent uniformly in $U$ because
\begin{align*}
&\|\lambda_2^{-m}\pi_2 F^{m}(x)-\lambda_2^{-n}\pi_2 F^n(x)\|
\\
&=
\|\{\lambda_2^{-m}\pi_2 F^{m}(x)-\lambda_2^{-n}\pi_2 F^n(x)\}
             -\{\lambda_2^{-m}\pi_2 F^{m}(O)-\lambda_2^{-n}\pi_2 F^n(O)\}\|
\\
&\le \sup_{\xi\in U}\|\lambda_2^{-m}D(\pi_2 F^{m})(\xi)-\lambda_2^{-n}D(\pi_2 F^n)(\xi)\|\, \|x\|,
\end{align*}
 for any $m\ne n\in\mathbb{N}$ by the Mean Value Theorem.
It means that the limit $
\psi:=\lim_{n\to\infty}\lambda_2^{-n}\pi_2 F^n $ gives a $C^1$
mapping with $D\psi(O)=(0,1)$.

Next, in order to prove that $\psi$ is also $C^{1,\alpha}$, we notice that
\begin{align}
\|D\psi(x)-D\psi(y)\|
&=\|\lim_{n\to \infty}\lambda_2^{-n}D(\pi_2 F^n)(x)-\lim_{n\to \infty}\lambda_2^{-n}D(\pi_2 F^n)(y)\|
\nonumber\\
&\le\sum_{i=1}^{\infty}\|\{\lambda_2^{-(n+1)}D(\pi_2 F^{n+1})(x)-\lambda_2^{-n}D(\pi_2 F^{n})(x)\}
\nonumber\\
 &\hspace{2.5cm}   -\{\lambda_2^{-(n+1)}D(\pi_2 F^{n+1})(y)-\lambda_2^{-n}D(\pi_2 F^{n})(y)\}\|
\nonumber\\
 &\hspace{1.3cm}+|\lambda_2|^{-1}\|D(\pi_2 F)(x)-D(\pi_2 F)(y)\|
 \label{Dpsi}
\end{align}
by (\ref{dfdf}). On the other hand, since
\begin{align*}
&p_1\cdots p_n-q_1\cdots q_n
\\
&=(p_1-q_1)p_2\cdots p_n+q_1(p_2-q_2)p_3\cdots p_n+\cdots +q_1\cdots q_{n-1}(p_n-q_n),
\end{align*}
we have
\begin{align}
&\|DF^n(x)-DF^n(y)\|
=
\Big\|\prod_{i=1}^{n}DF(F^{n-i}(x))-\prod_{i=1}^{n}DF(F^{n-i}(y))\Big\|
\nonumber\\
&\le\sum_{i=1}^{n}\bigg(\|DF(F^{n-1}(y))\|\cdots\|DF(F^{n-i+1}(y))\|
                          \,\|DF(F^{n-i}(x))-DF(F^{n-i}(y))\|
\nonumber\\
&~~~~~~~~~~\cdot\|DF(F^{n-i-1}(x))\|\cdots\|DF(x)\|\bigg)
\nonumber\\
&\le
\sum_{i=1}^{n}\bigg(\frac{(\lambda_2+\eta)^n}{\|DF(F^{n-i}(x))\|}
         \|DF(F^{n-i}(x))-DF(F^{n-i}(y))\|
\bigg)
\nonumber\\
&\le
\sum_{i=1}^{n}\bigg(\frac{(\lambda_2+\eta)^n}{\|DF(F^{n-i}(x))\|}
         L\sup_{\xi\in U}\|DF^{n-i}(\xi)\|^\alpha\,\|x-y\|^\alpha
\bigg)
\nonumber\\
&\le
\sum_{i=1}^{n}\bigg(\frac{(\lambda_2+\eta)^n}{\|DF(F^{n-i}(x))\|}
\, L(\lambda_2+\eta)^{(n-i)\alpha}\|x-y\|^\alpha\bigg)
\nonumber\\
&\le K_1(\lambda_2+\eta)^{n}\|x-y\|^\alpha.
\label{xany1ym}
\end{align}
It implies that
\begin{align}
& \|\{\lambda_2^{-(n+1)}D(\pi_2 F^{n+1})(x)-\lambda_2^{-n}D(\pi_2 F^{n})(x)\}
\nonumber\\
 &\hspace{2.5cm}   -\{\lambda_2^{-(n+1)}D(\pi_2 F^{n+1})(y)-\lambda_2^{-n}D(\pi_2 F^{n})(y)\}\|
\nonumber\\
& \le \|\lambda_2^{-(n+1)}\|\,\|\{DF^{n+1}(x)-\Lambda\, DF^{n}(x)\}
      -\{DF^{n+1}(y)-\Lambda\, DF^{n}(y)\}\|
\nonumber\\
& \le \|\lambda_2^{-(n+1)}\|\,\|\{DF(F^{n}(x))-\Lambda\}DF^{n}(x)
      -\{DF(F^{n}(y))-\Lambda\}DF^{n}(y)\}\|
\nonumber\\
& \le
\|\lambda_2^{-(n+1)}\|\,\|DF(F^n(x))-DF(F^n(y))\|\,\|DF^n(x)\|
\nonumber\\
&
\hspace{2.5cm}+\|\lambda_2^{-(n+1)}\|\,\|DF(F^n(y))-\Lambda\|\,\|DF^n(x)-DF^n(y)\|
\nonumber\\
& \le
K\bigg(\frac{(|\lambda_2|+\eta)^{1+\alpha}}{|\lambda_2|}\bigg)^n\|x-y\|^\alpha.
\label{bxzl2}
\end{align}
Hence, combining (\ref{Dpsi}) with (\ref{bxzl2}) we get
$
\|D\psi(x)-D\psi(y)\|\le L\|x-y\|^\alpha,
$
which proves that $\psi$ is $C^{1,\alpha}$.

Define a $C^{1,\alpha}$ diffeomorphism $\Theta_1: U\to \mathbb{R}^2$ by
$
\Theta_1(x):=(x_1, \psi(x)),
$
where $\psi$ is given above.
It is clear that $\Theta_1(O)=O$, $D\Theta_1(O)=$ id and its inverse $\Theta^{-1}$ satisfies $\psi\circ \Theta_1^{-1}(x)=x_2$.
Consequently, the $C^{1,\alpha}$ mapping $\widehat{F}:=\Theta_1\circ F\circ \Theta_1^{-1}$ satisfies
\begin{align*}
\pi_2 \widehat{F}(x)
&=\psi\circ F\circ \Theta_1^{-1}(x)=\lim_{n\to\infty}\lambda_2^{-n}\pi_2 F^{n+1}\circ \Theta_1^{-1}(x)
\\
&=\lambda_2\lim_{n\to\infty}\lambda_2^{-(n+1)}\pi_2 F^{n+1}\circ \Theta_1^{-1}(x)=
\lambda_2 \psi\circ \Theta_1^{-1}(x)=\lambda_2x_2.
\end{align*}
On the other hand, we know from \cite[Lemma 1]{WZWZ-JFA01} that
$\widehat{F}$ has a $C^{1,\alpha}$ invariant manifold
$
\Gamma:=\{(x_1,x_2)\in U: x_1=g(x_2)\},
$
where $g: U\cap \mathbb{R} \to \mathbb{R}$ is $C^{1,\alpha}$ such that
$g(0)=Dg(0)=0$. Thus, one can define a $C^{1,\alpha}$ diffeomorphism $\Theta_2:U\to\mathbb{R}^2$
by $\Theta_2(x):=(x_1-g(x_2), x_2)$, which satisfies that
$\Theta_2(O)=O$ and $D\Theta_2(O)=$ id. Then, putting
$\widetilde{F}:=\Theta_2\circ \widehat{F}\circ\Theta_2^{-1}$, we
obtain
\begin{eqnarray}
\left\{\begin{array}{ll}
\pi_1\widetilde{F}(x_1, x_2)=\pi_1\widehat{F}(x_1+g(x_2),x_2)-g(\pi_2\widehat{F}(x_1+g(x_2),x_2)),
\\
\pi_2\widetilde{F}(x_1, x_2)=\pi_2\widehat{F}(x_1+g(x_2),x_2)=\lambda_2x_2.
\end{array}\right.
\label{FwwCk}
\end{eqnarray}
Obviously,
\begin{eqnarray}
\pi_1\widetilde{F}(0,x_2)=\pi_1\widehat{F}(g(x_2),x_2)-g(\pi_2\widehat{F}(g(x_2),x_2))=0,
~~~~~ \forall x_2\in U\cap\mathbb{R},
\label{F1=0}
\end{eqnarray}
by the invariance of $\Gamma$. Hence, setting $\Theta:=\Theta_2\circ\Theta_1$ and combining (\ref{FwwCk}) with (\ref{F1=0}),
we complete the proof. \qquad$\Box$

The following are the main result of this section:

\begin{thm}
Let $F:\mathbb{R}^2\to \mathbb{R}^2$ be $C^{1,\alpha}$
with $\alpha\in (0,\alpha_0]$ such that {\rm (\ref{def-FR2})} and
{\rm (\ref{eig-val})} hold. Then there exists a homeomorphism $\Phi$
satisfies equation {\rm (\ref{schd-eqn})}, where $\Phi(O)=O$ and
$\Phi$ together with its inverse $\Phi^{-1}$ is differentiable at
$O$ such that $D\Phi(O)=D\Phi^{-1}(O)=$ {\rm id}.
\label{thm-DL}
\end{thm}

\vspace{-0.3cm}
{\bf Proof}.
According to Lemma \ref{lm-reduce}, we can consider the $C^{1,\alpha}$ mapping $\widetilde{F}:U \to \mathbb{R}^2$
defined in Lemma \ref{lm-reduce}, which satisfies (\ref{flatten}).
If the theorem holds for $\widetilde{F}$, then the theorem also holds for $F$ because the transformation $\Theta$ is $C^{1,\alpha}$ such that
$\Theta(O)=O$ and $D\Theta(O)=$ {\rm id}.

Notice that $D\widetilde{F}(O)=\Lambda={\rm diag}(\lambda_1,\lambda_2)$.
In view of (\ref{flatten}), it is reasonable to put
\begin{align}
D\widetilde{F}^n(x):=
\left(\begin{array}{cc}
a_n(x) & b_n(x)
\\
0 & \lambda_2^n
\end{array}\right),
\label{abcdn}
\end{align}
where $a_n:=\partial (\pi_1\widetilde{F}^n)/\partial x_1$ and $b_n:=\partial (\pi_1\widetilde{F}^n)/\partial x_2$,
which are both $C^{0,\alpha}$
functions such that $a_n(O)=\lambda_1^n$ and $b_n(O)=0$ for all $n\in \mathbb{N}$.
We claim that the sequence $\big(\lambda_1^{-n}a_n(x)\big)_{n\in\mathbb{N}}$ converges uniformly in $U$.
In fact, the equality
$$
D\widetilde{F}^n(x)=D\widetilde{F}(\widetilde{F}^{n-1}(x))D\widetilde{F}^{n-1}(x)
$$
gives
$
a_n(x) =a_1(\widetilde{F}^{n-1}(x))\,a_{n-1}(x)
$
and inductively
\begin{align*}
a_n(x)=\prod_{i=0}^{n-1}a_1(\widetilde{F}^i(x)), ~~~ \forall n\in\mathbb{N}.
\end{align*}
It follows that
\begin{align*}
|\lambda_1^{-(n+1)}a_{n+1}(x)-\lambda_1^{-n}a_n(x)|
&\le |\lambda_1|^{-(n+1)}|a_{n+1}(x)-\lambda_1a_n(x)|
\\
&\le |\lambda_1|^{-(n+1)}|a_1(\widetilde{F}^n(x))-\lambda_1|\prod_{i=0}^{n-1}|a_1(\widetilde{F}^i(x))|
\\
&\le K \bigg(\frac{(|\lambda_2|+\eta)^\alpha}{|\lambda_1|}\bigg)^n\,\prod_{i=0}^{n-1}\big(|\lambda_1|
       +L(|\lambda_2|+\eta)^{i\alpha}\big)
\\
&\le K (|\lambda_2|+\eta)^{n\alpha}\prod_{i=0}^{\infty}\big(1+(L/|\lambda_1|)(|\lambda_2|+\eta)^{i\alpha}\big)
\\
&\le M (|\lambda_2|+\eta)^{n\alpha},
\end{align*}
which implies the uniform convergence of the
sequence $\big(\lambda_1^{-n}a_n(x)\big)_{n\in\mathbb{N}}$
because $|\lambda_2|+\eta<1$ for small $\eta>0$. Therefore, the claim is proved.

The above claimed result enables us to assume that
\begin{align}
\lim_{n\to \infty} \lambda_1^{-n}a_n(x)= \varphi(x), ~~~~~ \forall x\in U,
\label{umfcon}
\end{align}
uniformly, where $\varphi:U\to \mathbb{R}$ is continuous such that $\varphi(O)=1$. According to the first equality of (\ref{flatten}), one sees that
$\pi_1\widetilde{F}^n(0,x_2)=0$ for all $n\in\mathbb{N}$ by induction and therefore
\begin{eqnarray}
\lim_{n\to\infty}\lambda_1^{-n}\pi_1\widetilde{F}^n(x)
=\lim_{n\to\infty}\int_{0}^{x_1}\lambda_1^{-n}{a}_n(t,x_2)dt
=\int_{0}^{x_1}\varphi(t,x_2)dt.
\label{limfn}
\end{eqnarray}
Define
\begin{align*}
\Phi(x):=\bigg(\lim_{n\to \infty}\lambda_1^{-n}\pi_1\widetilde{F}^n(x),x_2\bigg)
        =\bigg(\int_{0}^{x_1}\varphi(t,x_2)\,dt,\,x_2\bigg),
~~~~~ \forall x\in U.
\end{align*}
We can verify
\begin{align*}
\Phi\circ \widetilde{F}(x)&=\bigg(\lim_{n\to\infty}\lambda_1^{-n}\pi_1\widetilde{F}^{n+1}(x),\,\lambda_2x_2\bigg)
\\
&=\bigg(\lambda_1\lim_{n\to\infty}\lambda_1^{-(n+1)}\pi_1\widetilde{F}^{n+1}(x),\,\lambda_2x_2\bigg)=\Lambda\circ \Phi(x),
\end{align*}
implying that $\Phi$ is a solution of equation (\ref{schd-eqn}) such that
$\Phi(O)=O$.

In what follows, we prove that $\Phi$ is a homeomorphism and that $\Phi$ and its inverse $\Phi^{-1}$ are both
differentiable at $O$ such that $D\Phi(O)=D\Phi^{-1}(O)=$ {\rm id}.
The continuity of $\Phi$ is obvious.
Since the fact $\varphi(O)=1$ implies that $\varphi(x)>0$ for all $x\in U$,
we know that the function $\int_{0}^{x_1}\varphi(t,x_2)\,dt$
is strictly increasing with respect to the variable $x_1$.
Choose two different points $(x_1,x_2)$ and $(\tilde{x}_1,\tilde{x}_2)$ in $U$. Clearly,
$\Phi(x_1,x_2)\ne \Phi(\tilde{x}_1,\tilde{x}_2)$
by the definition if $x_2\ne \tilde{x}_2$.
If $x_2=\tilde{x}_2$ but $x_1\ne \tilde{x}_1$,
the strict monotonicity of $\int_{0}^{x_1}\varphi(t,x_2)\,dt$ with respect to $x_1$ also shows that
$\Phi(x_1,x_2)\ne \Phi(\tilde{x}_1,\tilde{x}_2)$.
Thus, $\Phi$ is invertible and therefore
$\Phi$ is a homeomorphism in $U$.
In order to prove the differentiability of $\Phi$ at $O$,
for an arbitrarily given $\varepsilon>0$, we notice that there
exists an integer $N(\varepsilon)>0$ such that
\begin{align*}
|\lambda_1^{-N(\varepsilon)}a_{N(\varepsilon)}(x)-\varphi(x)|\le \varepsilon/2, ~~~~~ \forall x\in U,
\end{align*}
because of the uniform convergence of the sequence
$(\lambda_1^{-n}a_n)_{n\in\mathbb{N}}$, given in (\ref{umfcon}).
Moreover, $|a_{N(\varepsilon)}(x)-a_{N(\varepsilon)}(y)|\le
L_{N(\varepsilon)}\|x-y\|^\alpha$ with a constant
$L_{N(\varepsilon)}>0$ since the function $a_n$ is $C^{0,\alpha}$
for every $n\in\mathbb{N}$. It follows that, for an arbitrarily
given $\varepsilon>0$, there is a
\begin{align*}
\delta(\varepsilon):=\big(\,\varepsilon|\lambda_1|^{N(\varepsilon)}/(2L_{N(\varepsilon)})\,\big)^{1/\alpha}>0
\end{align*}
such that if $\|x\|<\delta(\varepsilon)$ then
\begin{align*}
\frac{\|\Phi(x)-x\|}{\|x\|}&=\frac{|\int_{0}^{x_1}(\varphi(t,x_2)-1)\,dt|}{\|x\|}
\le |\varphi(\theta x_1,x_2)-1|
\\
&\le |\lambda_1^{-N(\varepsilon)}a_{N(\varepsilon)}(\theta x_1,x_2)-\varphi(\theta x_1,x_2)|+|\lambda_1^{-N(\varepsilon)}
a_{N(\varepsilon)}(\theta x_1,x_2)-1|
\\
&\le \varepsilon/2+|\lambda_1|^{-N(\varepsilon)}|a_{N(\varepsilon)}(\theta x_1,x_2)-\lambda_1^{N(\varepsilon)}|
\\
&= \varepsilon/2+|\lambda_1|^{-N(\varepsilon)}|a_{N(\varepsilon)}(\theta x_1,x_2)-a_{N(\varepsilon)}(O)|
\\
&\le \varepsilon/2+|\lambda_1|^{-N(\varepsilon)}L_{N(\varepsilon)}\|x\|^\alpha< \varepsilon,
\end{align*}
where $\theta\in (0,1)$.
% and $K(\varepsilon)>0$ is a number
%%
%\footnote{[ZWN130312] or `constant'?}
%%
%depending on $\varepsilon$
This implies that $\lim_{x\to
0}\|\Phi(x)-x\|/\|x\|=0$, i.e., $\Phi$ is differentiable at $O$ with
$D\Phi(O)=$ id.
On the other hand, setting
$y:=\Phi^{-1}(x)$, we have
\begin{align*}
\frac{\|\Phi^{-1}(x)-x\|}{\|x\|}
&=\frac{\|\Phi(y)-y\|}{\|y\|}\bigg(\frac{\|\Phi(y)\|}{\|y\|}\bigg)^{-1}\le\frac{\|\Phi(y)-y\|}{\|y\|}
     \bigg(\frac{\|y\|-\|\Phi(y)-y\|}{\|y\|}\bigg)^{-1}
\\
&\le \frac{\|\Phi(y)-y\|}{\|y\|}\bigg(1-\frac{\|\Phi(y)-y\|}{\|y\|}\bigg)^{-1}\to 0
\end{align*}
as $x\to 0$ (i.e., $y\to 0$), implying that $\Phi^{-1}$ is differentiable at $O$ with
$D\Phi^{-1}(O)=$ id. The proof is completed. \qquad$\Box$

\begin{re}
{\rm
$C^0$ linearization plus differentiability at the fixed point
was first investigated by van Strien (\cite{Stri-JDE90}) in 1990.
By solving a sequence of conjugacy equations, he proved that such a linearization can be realized for all $C^2$
hyperbolic diffeomorphisms in $\mathbb{R}^n$. Yet a problem of
the proof given in \cite{Stri-JDE90} were observed by Rayskin
(\cite{Ray-JDE98}). Therefore, van Strien's result was restated in
\cite{GuyHassRay-DCDS03} for $C^\infty$ hyperbolic diffeomorphisms
and proved via a different approach inspired from the normal form theory. Our Theorem \ref{thm-DL}
extends the previous results because we allow the smoothness of the
diffeomorphisms to be rather low, i.e., $C^{1,\alpha}$ for all
$\alpha\in (0,1]$. }
\end{re}

%%-----------------------------------------------------------------------------

\section{$C^{1,\beta}$ linearization in Siegel domain}
\setcounter{equation}{0}

This section is devoted to the case that $\lambda$ belongs to the Siegel domain, i.e., (\ref{sdcs}) holds.
As indicated in the Introduction, we need to decompose $F$ firstly by flattening it along its invariant
manifolds. Concretely, we transform $F$ into another mapping,
say $\widetilde{F}$, which together with its derivative is linearized on the invariant
manifolds (see (\ref{flat-1})-(\ref{flat-2}) below),
so that the two variables of $\widetilde{F}$ can be separated in some crucial estimates for the
sequences where iterates of $\widetilde{F}$ are involved. This decomposition will help us
eliminate the divergence caused by the part of expansion and get the convergence of the sequences.

Let $I$ and $J$, $I\subsetneq J\subset \mathbb{R}$,
be small compact intervals centered at $0$.
Since (\ref{FFbumpR2}) holds and (\ref{sdcs}) implies that the fixed point $O$ is hyperbolic,
by the Stable Manifold Theorem (cf. \cite{HPS-book77}) there exist a $C^{1,\alpha}$ stable manifold and an
unstable one of $F$.
Straightening up these manifolds, we may assume that
\begin{eqnarray}
\pi_{1}F(0,x_2)=0
\quad {\rm and} \quad
\pi_{2}F(x_1,0)=0,
~~~~~ \forall x_1,x_2\in \mathbb{R}.
\label{FFflat}
\end{eqnarray}
Then we get the following decomposition lemma.

\begin{lm}
Suppose that $F:\mathbb{R}^2\to \mathbb{R}^2$ is $C^{1,\alpha}$ such
that {\rm (\ref{def-FR2})}, {\rm (\ref{sdcs})} and
{\rm(\ref{FFflat})} hold. Then there exists a $C^{1,\alpha}$
diffeomorphism $\Psi:\mathbb{R}^2\to \mathbb{R}^2$ such that
$\Psi(O)=O$, $D\Psi(O)={\rm id}$, and
\begin{eqnarray}
(\Psi\circ F\circ \Psi^{-1})(x_1,0)=(\lambda_1x_1,0),
\qquad
D(\Psi\circ F\circ \Psi^{-1})(x_1,0)=\Lambda,
\label{flat-1}
\\
(\Psi\circ F\circ \Psi^{-1})(0,x_2)=(0,\lambda_2x_2),
\qquad
D(\Psi\circ F\circ \Psi^{-1})(0,x_2)=\Lambda~
\label{flat-2}
\end{eqnarray}
for all $x_1,x_2\in V\cap \mathbb{R}$.
\label{lm-flat}
\end{lm}

The proof of the lemma will be postponed to next section but we continue to give our main theorem of this
section as follows.

\begin{thm}
Suppose that $F:\mathbb{R}^2\to \mathbb{R}^2$ is $C^{1,\alpha}$ such
that {\rm (\ref{def-FR2})} and {\rm (\ref{sdcs})} hold. Then there
exists a $C^{1,\beta}$ diffeomorphism $\Phi: U\to \mathbb{R}^2$ such
that equation {\rm (\ref{schd-eqn})} holds, where
\begin{eqnarray*}
\beta=\left\{
\begin{array}{ll}
-{\sigma}(1-{\sigma})^{-1}\alpha, & ~~ {\it if} ~~ |\lambda_1\lambda_2|\le 1,
%~ ({\rm i.e.,} ~ \sigma\in[-1,0)),
\\
\vspace{-0.3cm}
\\
(1-{\sigma})^{-1}\alpha, & ~~ {\it if} ~~ |\lambda_1\lambda_2|> 1,
%~ ({\rm i.e.,} ~ \sigma\in(-\infty,-1))
\end{array}
\right.
\end{eqnarray*}
and $\sigma:={\log|\lambda_2|}/{\log|\lambda_1|}<0$.
\label{thm-C1-R2}
\end{thm}

{\bf Proof}. As mentioned above, one can always modify $F$ such that
(\ref{FFflat}) holds. Then, by Lemma \ref{lm-flat},
there exists a $C^{1,\alpha}$ diffeomorphism $\Psi:\mathbb{R}^2\to
\mathbb{R}^2$ such that $\Psi(O)=O$, $D\Psi(O)={\rm id}$ and
(\ref{flat-1})-(\ref{flat-2}) hold near $O$.
Multiplying the nonlinear part of $\Psi\circ
F\circ \Psi^{-1}$ by the bump function $\varrho$, we obtain a modified
$C^{1,\alpha}$ diffeomorphism $\widetilde{F}$ such that
$\widetilde{F}(O)=O$, $D\widetilde{F}(O)=\Lambda$ and
\begin{eqnarray}
%\widetilde{F}(x)=\Lambda x, \quad \forall x\in \mathbb{R}^2\backslash V,
%\qquad
\|D\widetilde{F}(x)-\Lambda\|\le \eta, \quad \forall x\in \mathbb{R}^2,
\label{FwwbumpR2}
\end{eqnarray}
where $\eta>0$ is a sufficiently small constant. Notice that $\widetilde{F}$ coincides with $\Psi\circ F\circ \Psi^{-1}$
in $U$. Furthermore, according to
(\ref{flat-1}), (\ref{flat-2}) and (\ref{FwwbumpR2}), one checks that
\begin{eqnarray}
\pi_1\widetilde{F}(0,x_2)=0,
\quad
D\widetilde{F}(0,x_2)=\Lambda,
\quad
\pi_2\widetilde{F}(x_1,0)=0,
\quad
D\widetilde{F}(x_1,0)=\Lambda
\label{flat-Fww}
\end{eqnarray}
for all $x_1,x_2\in \mathbb{R}$ and consequently
\begin{eqnarray}
D\widetilde{F}^{-1}(0,x_2)=\Lambda^{-1},
\qquad
D\widetilde{F}^{-1}(x_1,0)=\Lambda^{-1},
~~~~~ \forall x_1,x_2\in \mathbb{R}.
\label{flat-Fww-1}
\end{eqnarray}

Next, we claim the uniform convergence of the sequences
\begin{eqnarray*}
\Big(\lambda_1^nD(\pi_1\widetilde{F}^{-n})(x)\Big)_{n\in\mathbb{N}}
~~~~~{\rm and}~~~~~
\Big(\lambda_2^{-n}D(\pi_2\widetilde{F}^n)(x)\Big)_{n\in\mathbb{N}}
\end{eqnarray*}
in $U$. If the claim is true, then
$\widetilde{\Phi}:=(\varphi_1,\varphi_2)$ is $C^1$ such that
$\widetilde{\Phi}\circ \widetilde{F}=\Lambda\circ \widetilde{\Phi}$, where
\begin{eqnarray}
\varphi_1:=\lim_{n\to\infty}\lambda_1^{n}\pi_1\widetilde{F}^{-n}
\quad{\rm and}\quad
\varphi_2:=\lim_{n\to\infty}\lambda_2^{-n}\pi_2\widetilde{F}^n.
\label{phi12}
\end{eqnarray}
Indeed,
\begin{eqnarray*}
\widetilde{\Phi}\circ \widetilde{F}
=\left(
\begin{array}{ll}
\lim_{n\to\infty}\lambda_1^{n}\pi_1\widetilde{F}^{-n}\circ \widetilde{F}
\\
\lim_{n\to\infty}\lambda_2^{-n}\pi_2\widetilde{F}^n\circ \widetilde{F}
\end{array}
\right)
=
\left(
\begin{array}{ll}
\lambda_1\lim_{n\to\infty}\lambda_1^{n-1}\pi_1\widetilde{F}^{-(n-1)}
\\
\lambda_2\lim_{n\to\infty}\lambda_2^{-(n+1)}\pi_2\widetilde{F}^{n+1}
\end{array}
\right)
=\Lambda \circ\widetilde{\Phi}.
\end{eqnarray*}
 Moreover we have
$
D\widetilde{\Phi}(O)={\rm id},
$
which implies that
$\widetilde{\Phi}$ is a $C^1$ diffeomorphism near $O$. Therefore,
$F$ can be linearized by the $C^1$ diffeomorphism $\Phi:=\widetilde{\Phi}\circ \Psi$ near $O$.

In order to show the uniform convergence of the two sequences given in (\ref{phi12}), we claim that
\begin{eqnarray}
|\pi_1\widetilde{F}^n(x)|\le M |\lambda_1|^n|x_1|,
~~~ |\pi_2\widetilde{F}^{-n}(x)|\le M |\lambda_2|^{-n}|x_2|,
~~~~ \forall n\in\mathbb{N}, \,\forall x\in U.
\label{F1F1'}
\end{eqnarray}
In fact, (\ref{FwwbumpR2}) and (\ref{flat-Fww}) give
\begin{eqnarray*}
|\pi_1 \widetilde{F}(x)|
&=&|\pi_1 \widetilde{F}(x_1,x_2)-\pi_1\widetilde{F}(0,x_2)|
    \le\sup_{\xi\in (0,x_1)}|\partial_{x_1}(\pi_1\widetilde{F})(\xi,x_2)|\,|x_1|
\nonumber\\
&\le&\Big(|\lambda_1|+\sup_{\xi\in (0,x_1)}|\partial_{x_1}(\pi_1\widetilde{F})(\xi,x_2)-\lambda_1|\Big)|x_1|
\nonumber\\
&\le& (|\lambda_1|+\eta)|x_1|, ~~~~~ \forall x\in \mathbb{R}^2,
%\label{FSF1}
\end{eqnarray*}
for an $\eta<(1-|\lambda_1|)/2$. It follows that
\begin{eqnarray}
|\pi_1 \widetilde{F}^i(x)|\le (|\lambda_1|+\eta)^i|x_1|\le
(|\lambda_1|+\eta)^i, ~~~~~ \forall i\in \mathbb{N},~ \forall x\in U,
\label{Fn}
\end{eqnarray}
by induction. On the other hand, we also have
\begin{eqnarray}
|\pi_1 \widetilde{F}(x)|
&\le&\Big(|\lambda_1|+\sup_{\xi\in (0,x_1)}|\partial_{x_1}(\pi_1\widetilde{F})(\xi,x_2)-\lambda_1|\Big)|x_1|
\nonumber\\
&=&\Big(|\lambda_1|+\sup_{\xi\in (0,x_1)}|\partial_{x_1}(\pi_1\widetilde{F})(\xi,x_2)-\partial_{x_1}(\pi_1\widetilde{F})(0,x_2)|\Big)|x_1|
\nonumber\\
&\le& (|\lambda_1|+L|x_1|^\alpha)|x_1|
\label{FSF2}
\end{eqnarray}
since $D\widetilde{F}$ is $C^{0,\alpha}$. Thus, substituting
(\ref{Fn}) in (\ref{FSF2}), we get
\begin{eqnarray*}
|\pi_1 \widetilde{F}^{i+1}(x)|\le (|\lambda_1|+L (|\lambda_1|+\eta)^{i{\alpha}})|\pi_1\widetilde{F}^{i}(x)|,
~~~ \forall i\in\mathbb{N}\cup \{0\}, ~ \forall x\in U,
\end{eqnarray*}
which proves the first inequality of (\ref{F1F1'})
by induction because
\begin{eqnarray}
\prod_{i=0}^{\infty}\Big\{1+\frac{L}{|\lambda_1|}(|\lambda_1|+\eta)^{i{\alpha}}\Big\}<\infty.
\label{smyd}
\end{eqnarray}
Similarly, we can also prove the second inequality of (\ref{F1F1'}). Moreover, one sees that
\begin{eqnarray}
\|D\widetilde{F}^n(x)\|\le M |\lambda_2|^n,
~~~
\|D\widetilde{F}^{-n}(x)\|\le M |\lambda_1|^{-n},
~~~~ \forall n\in \mathbb{N},\, \forall x\in U.
\label{DFn}
\end{eqnarray}
The first inequality holds because, by
(\ref{F1F1'}),
\begin{align*}
\|D\widetilde{F}^n(x)\|
&\le \prod_{i=0}^{n-1}\|D\widetilde{F}(\widetilde{F}^i(x))\|
\le\prod_{i=0}^{n-1}\{\|\Lambda\|+\|D\widetilde{F}(\widetilde{F}^i(x))-D\widetilde{F}(0,\pi_2\widetilde{F}^i(x))\|\}
\\
&\le\prod_{i=0}^{n-1}(|\lambda_2|+L|\pi_1\widetilde{F}^i(x)|^\alpha)
     \le\prod_{i=0}^{n-1}(|\lambda_2|+LM^\alpha|\lambda_1|^{i\alpha})
\end{align*}
and
$\prod_{i=0}^{\infty}\{1+LM^\alpha|\lambda_1|^{i\alpha}/|\lambda_2|\}<\infty$.
The second one holds for the same reason. Then, in view of (\ref{flat-Fww-1}),
(\ref{F1F1'}) and (\ref{DFn}), we get
\begin{eqnarray*}
&&\|\lambda_1^{n+1}D(\pi_{1}\widetilde{F}^{-(n+1)})(x)-\lambda_1^{n}D(\pi_{1}\widetilde{F}^{-n})(x)\|
\\
&&\le|\lambda_1^{n+1}|\,\|D\widetilde{F}^{-1}(\widetilde{F}^{-n}(x))-\Lambda^{-1}\|\,\|D\widetilde{F}^{-n}(x)\|
\\
&&\le|\lambda_1^{n+1}|\,\|D\widetilde{F}^{-1}(\widetilde{F}^{-n}(x))
               -D\widetilde{F}^{-1}(\pi_1\widetilde{F}^{-n}(x),0)\|\,\|D\widetilde{F}^{-n}(x)\|
\\
&&\le L|\lambda_1^{n+1}|\,|\pi_{2}\widetilde{F}^{-n}(x)|^\alpha\|D\widetilde{F}^{-n}(x)\|
\\
&&\le K|\lambda_2|^{-n\alpha}.
\end{eqnarray*}
This implies the uniform convergence of the sequence
$(\lambda_1^nD(\pi_1\widetilde{F}^{-n})(x))_{n\in\mathbb{N}}$ by an analogue of equality (\ref{dfdf}) since $|\lambda_2|^{-\alpha}<1$.
Also we have
\begin{eqnarray*}
&&\|\lambda_2^{-(n+1)}D(\pi_{2}\widetilde{F}^{n+1})(x)-\lambda_2^{-n}D(\pi_{2}\widetilde{F}^{n})(x)\|
\\
&&\le|\lambda_2^{-(n+1)}|\,\|D\widetilde{F}(\widetilde{F}^{n}(x))-\Lambda\|\,\|D\widetilde{F}^{n}(x)\|
\\
&&\le|\lambda_2^{-(n+1)}|\,\|D\widetilde{F}(\widetilde{F}^{n}(x))
         -D\widetilde{F}(0,\pi_2\widetilde{F}^{n}(x))\|\,\|D\widetilde{F}^{n}(x)\|
\\
&&\le L|\lambda_2^{-(n+1)}|\,|\pi_{1}\widetilde{F}^{n}(x)|^\alpha\|D\widetilde{F}^{n}(x)\|
\\
&&\le K|\lambda_1|^{n\alpha},
\end{eqnarray*}
implying the uniform convergence of $(\lambda_2^{-n}D(\pi_2\widetilde{F}^n)(x))_{n\in\mathbb{N}}$.
Hence $F$ can be linearized by a $C^1$ diffeomorphism $\Phi=\widetilde{\Phi}\circ \Psi$ near $O$,
as indicated in the previous paragraph.

Next, we show that the diffeomorphism $\Phi=\widetilde{\Phi}\circ \Psi$ is not only $C^1$
but also $C^{1,\beta}$ near $O$ for some $\beta\in (0,\alpha]$.
Since Lemma \ref{lm-flat} implies that $\Psi$ is $C^{1,\alpha}$, it suffices
to prove that $\widetilde{\Phi}=(\varphi_1(x),\varphi_2(x))$ is $C^{1,\beta}$.
For this purpose, we first investigate the function $\varphi_1=\lim_{n\to\infty}\lambda_1^{n}\pi_1\widetilde{F}^{-n}$. Similarly to
(\ref{Dpsi}) and (\ref{bxzl2}), one computes that
\begin{align}
&\|D\varphi_1(x)-D\varphi_1(y)\|
    =\|\lim_{n\to\infty}\lambda_1^{n}D(\pi_1\widetilde{F}^{-n})(x)-\lim_{n\to\infty}\lambda_1^{n}D(\pi_1\widetilde{F}^{-n})(y)\|
\nonumber\\
&\le \sum_{n=1}^{\infty}|\lambda_1^{n+1}|\Big(\|D\widetilde{F}^{-1}(\widetilde{F}^{-n}(x))
           -D\widetilde{F}^{-1}(\widetilde{F}^{-n}(y))\|\,\|D\widetilde{F}^{-n}(x)\|
\nonumber\\
&~~~~~   +\|D\widetilde{F}^{-1}(\widetilde{F}^{-n}(y))-\Lambda^{-1}\|\,\|D\widetilde{F}^{-n}(x)-D\widetilde{F}^{-n}(y)\|\Big)
        +L|\lambda_1|\,\|x-y\|^\alpha.
\label{hold-est}
\end{align}
In order to estimate the term
\begin{eqnarray*}
\sum_{n=1}^{\infty}|\lambda_1^{n+1}|\,\|D\widetilde{F}^{-1}(\widetilde{F}^{-n}(x))
     -D\widetilde{F}^{-1}(\widetilde{F}^{-n}(y))\|\,\|D\widetilde{F}^{-n}(x)\|
\end{eqnarray*}
given in (\ref{hold-est}), we observe that either
\begin{eqnarray}
&&\|D\widetilde{F}^{-1}(\widetilde{F}^{-n}(x))-D\widetilde{F}^{-1}(\widetilde{F}^{-n}(y))\|
\nonumber\\
&&\le L\|\widetilde{F}^{-n}(x)-\widetilde{F}^{-n}(y)\|^\alpha
    \le L\sup_{\xi\in U}\|D\widetilde{F}^{-n}(\xi)\|^\alpha\|x-y\|^{\alpha}
\nonumber\\
&&\le K|\lambda_1|^{-n\alpha}\|x-y\|^{\alpha}
\label{anxany1}
\end{eqnarray}
by (\ref{DFn}) or
\begin{eqnarray}
&&\|D\widetilde{F}^{-1}(\widetilde{F}^{-n}(x))-D\widetilde{F}^{-1}(\widetilde{F}^{-n}(y))\|
\nonumber\\
&&\le \|D\widetilde{F}^{-1}(\widetilde{F}^{-n}(x))-\Lambda^{-1}\|+\|D\widetilde{F}^{-1}(\widetilde{F}^{-n}(y))-\Lambda^{-1}\|
\nonumber\\
&&\le L (|\pi_{2}\widetilde{F}^{-n}(x)|^\alpha+|\pi_{2}\widetilde{F}^{-n}(y)|^\alpha)\le K|\lambda_2|^{-n\alpha}
\label{anxany2}
\end{eqnarray}
by (\ref{F1F1'}) for all $x,y\in U$. Choose
$n_1(x,y):=\log \|x-y\|/\log (|\lambda_1|/|\lambda_2|)$. Obviously
$n_1(x,y)$, which is simply denoted by $n_1$, is larger than $1$ for small $x$ and $y$ and satisfies that
$%\begin{eqnarray*}
|\lambda_1|^{-n_1}\|x-y\|=|\lambda_2|^{-n_1}.
$%\end{eqnarray*}
Therefore,
\begin{eqnarray}
|\lambda_1|^{n_1}=\|x-y\|^{(1-{\sigma})^{-1}}
\quad \mbox{and}\quad
|\lambda_2|^{n_1}=\|x-y\|^{{\sigma}(1-{\sigma})^{-1}},
\label{1212}
\end{eqnarray}
where $\sigma$, defined as in the statement of Theorem~\ref{thm-C1-R2}, satisfies that $|\lambda_2|=|\lambda_1|^\sigma$.
Then, (\ref{anxany1})-(\ref{1212}) yield
\begin{eqnarray}
&&\sum_{n=1}^{\infty}|\lambda_1^{n+1}|\,\|D\widetilde{F}^{-1}(\widetilde{F}^{-n}(x))
          -D\widetilde{F}^{-1}(\widetilde{F}^{-n}(y))\|\,\|D\widetilde{F}^{-n}(x)\|
\nonumber\\
&&\le |\lambda_1|MK\sum_{n=1}^{[n_1]}|\lambda_1|^{-n\alpha}\|x-y\|^{\alpha}
          +|\lambda_1|MK\sum_{n=[n_1]+1}^{\infty}|\lambda_2|^{-n\alpha}
\nonumber\\
&&\le M_1\|x-y\|^{-{\sigma}(1-{\sigma})^{-1}\alpha}, ~~~~~ \forall x,y\in U,
\label{yymxwm}
\end{eqnarray}
where $[n_1]$ denotes the greatest integer not exceeding $n_1$.
On the other hand, in order to estimate the term
\begin{eqnarray*}
\sum_{n=1}^{\infty}|\lambda_1^{n+1}|\,\|D\widetilde{F}^{-1}(\widetilde{F}^{-n}(y))-\Lambda^{-1}\|\,
         \|D\widetilde{F}^{-n}(x)-D\widetilde{F}^{-n}(y)\|
\end{eqnarray*}
given in (\ref{hold-est}), we observe that either
\begin{align}
&\|D\widetilde{F}^{-n}(x)-D\widetilde{F}^{-n}(y)\|
\nonumber\\
&\le
\sum_{i=1}^{n}\bigg(\frac{M|\lambda_1|^{-n}}{\|D\widetilde{F}^{-1}(\widetilde{F}^{-(n-i)}(x))\|}
\,
LM^{\alpha}|\lambda_1|^{-(n-i){\alpha}}\|x-y\|^{{\alpha}}\bigg)
\nonumber\\
&\le K|\lambda_1|^{-n(1+\alpha)}\|x-y\|^{{\alpha}}
\label{xany1}
\end{align}
(the calculation in (\ref{xany1}) is similar to (\ref{xany1ym})) or
\begin{eqnarray}
\|D\widetilde{F}^{-n}(x)-D\widetilde{F}^{-n}(y)\|
\le \|D\widetilde{F}^{-n}(x)\|+\|D\widetilde{F}^{-n}(y)\|\le K|\lambda_1|^{-n}
\label{xany2}
\end{eqnarray}
for all $x,y\in U$. Choose
$n_2(x,y):=\log \|x-y\|/\log (|\lambda_1|^{1+\alpha}/|\lambda_2|)$. Obviously $n_2(x,y)$, which is
simply denoted by $n_2$, is larger than $1$ for small $x$ and $y$ and satisfies that
$%\begin{eqnarray*}
|\lambda_1|^{-(1+\alpha)n_2}\|x-y\|^{{\alpha}}=|\lambda_1|^{-n_2}.
$%\end{eqnarray*}
Therefore,
\begin{eqnarray}
|\lambda_1|^{n_2}=\|x-y\|
\quad\mbox{and}\quad
|\lambda_2|^{n_2}=\|x-y\|^{\sigma}.
\label{yyy}
\end{eqnarray}
Then, (\ref{xany1})-(\ref{yyy}) yield
\begin{eqnarray}
&&\sum_{n=1}^{\infty}|\lambda_1^{n+1}|\,\|D\widetilde{F}^{-1}(\widetilde{F}^{-n}(y))-\Lambda^{-1}\|
               \,\|D\widetilde{F}^{-n}(x)-D\widetilde{F}^{-n}(y)\|
\nonumber\\
&&\le|\lambda_1|KLM^\alpha\sum_{n=1}^{[n_2]}|\lambda_1\lambda_2|^{-n\alpha}\|x-y\|^{{\alpha}}
          +|\lambda_1|KLM^\alpha\sum_{n=[n_2]+1}^{\infty}|\lambda_2|^{-n\alpha}
\nonumber\\
&&\le M_2\|x-y\|^{\varsigma_1},~~~~~ \forall x,y\in U,
\label{wngyo}
\end{eqnarray}
where $\varsigma_1:=\min\{-{\sigma}\alpha,~\alpha-\gamma\}$ and
$\gamma>0$ is a sufficiently small constant.
Note that in the last row of formula (\ref{wngyo}) we need to use the inequality
\begin{eqnarray*}
\sum_{n=1}^{[n_2]}|\lambda_1\lambda_2|^{-n\alpha}
\le\left\{
\begin{array}{ll}
K_1, &~~ {\rm if} ~ |\lambda_1\lambda_2|>1,
\\
K_2 \|x-y\|^{-\gamma}, &~~ {\rm if} ~ |\lambda_1\lambda_2|=1,
\\
K_3|\lambda_1\lambda_2|^{-n_2\alpha}=K_3\|x-y\|^{-(1+\sigma)\alpha}, &~~ {\rm if} ~ |\lambda_1\lambda_2|< 1,
\end{array}
\right.~~~
\end{eqnarray*}
which is true becuase
\begin{align*}
&\sum_{n=1}^{[n_2]}|\lambda_1\lambda_2|^{-n\alpha}
=
\frac{|\lambda_1\lambda_2|^{-\alpha}}{1-|\lambda_1\lambda_2|^{-\alpha}}
   (1-|\lambda_1\lambda_2|^{-[n_2]\alpha})
\\
&~~\le\left\{
\begin{array}{ll}
|\lambda_1\lambda_2|^{-\alpha}(1-|\lambda_1\lambda_2|^{-\alpha})^{-1}=:K_1,
&~~ {\rm if} ~ |\lambda_1\lambda_2|>1,
\\
\vspace{-0.3cm}
\\
|\lambda_1\lambda_2|^{-\alpha}(|\lambda_1\lambda_2|^{-\alpha}-1)^{-1}|\lambda_1\lambda_2|^{-[n_2]\alpha}
     =:K_3|\lambda_1\lambda_2|^{-[n_2]\alpha},
     &~~ {\rm if} ~ |\lambda_1\lambda_2|< 1,
\end{array}
\right.
\end{align*}
when $|\lambda_1\lambda_2|\ne 1$ and
\begin{align*}
\sum_{n=1}^{[n_2]}|\lambda_1\lambda_2|^{-n\alpha}
&=[n_2]\le n_2=\log \|x-y\|/\log (|\lambda_1|^{1+\alpha}/|\lambda_2|)
\\
&\le -K_2\log\|x-y\|=K_2\log\|x-y\|^{-1}\le K_2 \|x-y\|^{-\gamma}
~~~~~~~~~~~~~
\end{align*}
when $|\lambda_1\lambda_2|=1$
since $\lim_{\|x-y\|\to 0}\|x-y\|^\gamma\log \|x-y\|^{-1}=0$ for $\gamma>0$.
Combining (\ref{hold-est}), (\ref{yymxwm}) with (\ref{wngyo}), we get
\begin{eqnarray}
|\varphi_1(x)-\varphi_1(y)|
\le L_1\|x-y\|^{\varsigma_2}
%\nonumber\\
%&=&
=L_1\|x-y\|^{-{\sigma}(1-{\sigma})^{-1}\alpha},
~~ \forall x,y\in U,
\label{betaes1}
\end{eqnarray}
where $\varsigma_2:=\min\{-{\sigma}(1-{\sigma})^{-1}\alpha,~-{\sigma}\alpha,~\alpha-\gamma\}$.

For the function $\varphi_2$ given in (\ref{phi12}), replacing $\widetilde{F}$, $\lambda_1$, $\lambda_2$, $\pi_1$ and $\sigma$
with $\widetilde{F}^{-1}$, $\lambda_2^{-1}$, $\lambda_1^{-1}$, $\pi_2$ and $\sigma^{-1}$ respectively in the above proof for $\varphi_1$,
we get
\begin{eqnarray}
|\varphi_2(x)-\varphi_2(y)|\le L_2\|x-y\|^{(1-{\sigma})^{-1}\alpha},
\qquad \forall x,y\in U.
\label{betaes2}
\end{eqnarray}
Thus, from (\ref{betaes1}) and (\ref{betaes2}) we have
$$
\|\Phi(x)-\Phi(y)\|\le L \|x-y\|^{\varsigma_3},
$$
where $\varsigma_3:=\min\{-\sigma(1-{\sigma})^{-1}\alpha,~(1-{\sigma})^{-1}\alpha\}$,
and the proof is completed by putting $\beta:=\min\{-\sigma(1-{\sigma})^{-1}\alpha,(1-{\sigma})^{-1}\alpha\}$.
\qquad$\Box$

\begin{re}
{\rm
In the above proof an analogous strategy to \cite{WZWZ-JFA01} is employed when we dealt with planar contractions,
but it needs some essential improvements. This fact
can be observed by a simple example
that the infinite product
$
\prod_{i=1}^{\infty}(1+L(|\lambda_2|+\eta)^{i\alpha})
$
with a small constant $\eta>0$ converges when
$|\lambda_2|<1$ but diverges when $|\lambda_2|>1$, where the
appearance of $(|\lambda_2|+\eta)$ in the infinite product
is due to the inequality
\begin{eqnarray*}
\|D(\pi_1F)(F^i(x))-D(\pi_1F)(O)\|\le L\|F^i(x)\|^\alpha\le
L(|\lambda_2|+\eta)^{i\alpha}.
\end{eqnarray*}
Thus we are forced to flatten $F$ along the invariant manifolds for a decomposition so that the
above given constant $|\lambda_2|$ ($>1$) can be replaced with $|\lambda_1|$ ($<1$), as seen in (\ref{smyd}).
}
\end{re}

\begin{re}
{\rm
It is hard to strengthen the $C^1$ linearization to $C^{1,\beta}$ linearization by using
the method of \cite{AraBeliZhu-book96}.
In fact, in
\cite{AraBeliZhu-book96} equation (\ref{schd-eqn}) is decomposed
into three functional equations, two of which depend on only the
first variable and only the second variable separately and the other
of which is equal to 0 on both axes. In order to solve the third
equation, a special norm
\begin{eqnarray*}
\|h\|_{\alpha}:=\max\Big\{\sup_{x_2\ne 0}\frac{|\partial_{x_1}h(x_1,x_2)|}{|x_2|^\alpha},~
\sup_{x_1\ne 0}\frac{|\partial_{x_2}h(x_1,x_2)|}{|x_1|^\alpha}\Big\}
\end{eqnarray*}
with $\alpha\in(0,1]$ is employed so as to obtain a contraction constant $\iota\in (0,1)$ in application of the
well-known Banach's contraction principle. However, such a constant $\iota$ cannot be deduced
if the norm $\|\cdot\|_\alpha$ is changed into the norm
$\|\cdot\|_{C^{1,\beta}}$ defined in (\ref{Cnorm}) with $\alpha:=\beta$.
This makes their method not applicable.
}
\end{re}

\begin{re}
{\rm
One can also tackle the linearization problem by means of
invariant foliations of mappings (cf. \cite{Tan-JDE00}). However,
such a method is not available for $C^1$ linearization of
$C^{1,\alpha}$ mappings with $\alpha\in (0,1)$. One of the main
reasons can be observed from the Lyapunov-Perron equation (cf. \cite[Lemma 3.3]{ChHaTan-JDE97})
\begin{align}
v_n(x,y_1)=~&\lambda_1^n(y_1-\pi_1
x)+\sum_{k=0}^{n-1}\lambda_1^{n-k-1}\{\pi_1 R(v_k+F^k(x))-\pi_1
R(F^k(x))\}
\nonumber\\
&-\sum_{k=n}^{\infty}\lambda_2^{n-k-1}\{\pi_2 R(v_k+F^k(x))-\pi_2
R(F^k(x))\}, ~~~ \forall n\ge 0,
\label{eqns-foli}
\end{align}
which plays a key role in proving the existence of invariant foliations. Here $v_n$'s are unknown and $R:=F-\Lambda$. Indeed,
an operator $T$ defined by derivatives of the right hand side of (\ref{eqns-foli}) can be proved to satisfy an inequality
\begin{align*}
\|T((v_k)_{k\in\mathbb{N}})-T((\tilde{v}_k)_{k\in\mathbb{N}})\|<\theta\|(v_k)_{k\in\mathbb{N}}-(\tilde{v}_k)_{k\in\mathbb{N}}\|^\alpha
\end{align*}
with $\theta\in (0,1)$ since $DR$ is $C^{0,\alpha}$. Then the Banach's contraction principle fails to solve equations
(\ref{eqns-foli}) because $\alpha\in (0,1)$.
}
\end{re}

%%------------------------------------------------

\section{Decomposition of the mapping $F$}
\setcounter{equation}{0}

According to the above proof of Theorem \ref{thm-C1-R2}, one sees that Lemma \ref{lm-flat} plays an important role
as it is actually a decomposition of $F$ along its invariant manifolds. The decomposition helps us
prove the convergence of the sequences defined in (\ref{phi12}).
In this section, we complementarily give the proof
of our Lemma \ref{lm-flat}.

First of all, we take into account the two equations given in (\ref{flat-1}).
Observe that if we can find a mapping $\Psi$ fixing $O$ such that the second equation of
(\ref{flat-1}) holds then $\Psi$ automatically satisfies
the first equation of (\ref{flat-1}) because
\begin{eqnarray*}
\pi_1(\Psi\circ F\circ \Psi^{-1})(x_1,0)=\int_{0}^{x_1}\partial_{x_1}(\pi_1\Psi\circ F\circ \Psi^{-1})(\tau,0)d\tau=\lambda_1x_1
\end{eqnarray*}
and
\begin{eqnarray*}
\pi_2(\Psi\circ F\circ \Psi^{-1})(x_1,0)=\int_{0}^{x_1}\partial_{x_1}(\pi_2\Psi\circ F\circ \Psi^{-1})(\tau,0)d\tau=0.
\end{eqnarray*}
Concerning the second equation of (\ref{flat-1}), it suffices to discuss the equation
\begin{eqnarray}
D\Psi(F(x_1,0))\,DF(x_1,0)=\Lambda\, D\Psi(x_1,0)
\label{eqn-1}
\end{eqnarray}
and find a local $C^{1,\alpha}$ solution $\Psi$ of (\ref{eqn-1}) such that
\begin{eqnarray}
\Psi(O)=O, ~~~ D\Psi(O)={\rm id}, ~~~ \pi_{2}\Psi(x_1,0)=0.
\label{PFPflat}
\end{eqnarray}
In fact, (\ref{PFPflat}) gives
$\Psi^{-1}(O)=O$, $D\Psi^{-1}(O)={\rm id}$ and $\pi_{2}\Psi^{-1}(x_1,0)=0$. This makes it reasonable to
replace $(x_1,0)$ with $\Psi^{-1}(x_1,0)$ in (\ref{eqn-1}) and get
\begin{eqnarray*}
D\Psi(F\circ \Psi^{-1}(x_1,0))\,DF(\Psi^{-1}(x_1,0))\,\{D\Psi(\Psi^{-1}(x_1,0))\}^{-1}=\Lambda.
\end{eqnarray*}
On the other hand, note that
\begin{align*}
&D\Psi(F\circ \Psi^{-1}(x_1,0))\,DF(\Psi^{-1}(x_1,0))\,\{D\Psi(\Psi^{-1}(x_1,0))\}^{-1}
\\
&=D\Psi(F\circ \Psi^{-1}(x_1,0))\,DF(\Psi^{-1}(x_1,0))\,D\Psi^{-1}(x_1,0)
    =D(\Psi\circ F\circ \Psi^{-1})(x_1,0).
\end{align*}
Hence we have $D(\Psi\circ F\circ \Psi^{-1})(x_1,0)=\Lambda$, i.e., $\Psi$ also satisfies
the second equation of (\ref{flat-1}) for small $x_1$.

In order to solve equation (\ref{eqn-1}), in view of (\ref{FFflat}) and (\ref{PFPflat}), we put
\begin{eqnarray}
DF(x_1,0)=
\left(\begin{array}{cc}
a_{11}(x_1) & a_{12}(x_1)
\\
0 & a_{22}(x_1)
\end{array}\right),
\quad
D\Psi(x_1,0)=
\left(\begin{array}{cc}
p_{11}(x_1) & p_{12}(x_1)
\\
0 & p_{22}(x_1)
\end{array}\right)
\label{FP-matrix}
\end{eqnarray}
for all $x_1\in\mathbb{R}$, where $a_{ij}$'s and $p_{ij}$'s are functions defined on $\mathbb{R}$.
Let $f:\mathbb{R}\to\mathbb{R}$ be defined by
\begin{eqnarray*}
f(s)=\pi_{1}F(s,0),
\qquad \forall s\in \mathbb{R}.
\end{eqnarray*}
In view of (\ref{def-FR2}), (\ref{FFbumpR2}) and (\ref{FP-matrix}), the above given diffeomorphism $f:\mathbb{R}\to\mathbb{R}$ is $C^{1,\alpha}$
such that
\begin{eqnarray}
&&f(0)=(0),
\qquad
Df(0)=\lambda_1,
\qquad
Df(s)=a_{11}(s),
\label{f-cond}
\\
&&f(s)=\lambda_1 s, \quad \forall s\in \mathbb{R}\backslash J,
\quad {\rm and} \quad
|Df(s)-\lambda_1|\le \eta, \quad \forall s\in \mathbb{R},
\label{fsbump}
\end{eqnarray}
for a sufficiently small $\eta>0$ depending on $J$.
Then equation (\ref{eqn-1}) can be decomposed into the following three functional equations:
\begin{eqnarray}
&&a_{11}(s)p_{11}(f(s))=\lambda_1p_{11}(s),
\label{FE-1}
\\
&&a_{22}(s)p_{12}(f(s))=\lambda_1p_{12}(s)-a_{12}(s)p_{11}(f(s)),
\label{FE-2}
\\
&&a_{22}(s)p_{22}(f(s))=\lambda_2p_{22}(s),
\label{FE-3}
\end{eqnarray}
where $f$ and $a_{ij}$'s are given and $p_{ij}$'s are unknown.
For equations (\ref{FE-1})-(\ref{FE-3}) we have the following lemma:
\begin{lm}
Let $f:\mathbb{R}\to \mathbb{R}$ be a $C^{1,\alpha}$ diffeomorphism satisfying {\rm (\ref{f-cond})}-{\rm (\ref{fsbump})}
and let $a_{ij}$'s be $C^{0,\alpha}$ functions given in {\rm (\ref{FP-matrix})}.
Then equations {\rm (\ref{FE-1})-(\ref{FE-3})} have $C^{0,\alpha}$ solutions $p_{ij}$'s defined on $I$ such that
\begin{eqnarray}
p_{11}(0)=p_{22}(0)=1
\quad {\rm and} \quad
p_{12}(0)=0.
\label{p1234}
\end{eqnarray}
\label{lm-3FE}
\end{lm}

{\bf Proof}. For convenience, we rewrite equations {\rm (\ref{FE-1})-(\ref{FE-3})} as
\begin{eqnarray}
&&\tilde{p}_{11}(s)=\frac{a_{11}(s)}{\lambda_1}\tilde{p}_{11}(f(s))+\frac{a_{11}(s)}{\lambda_1}-1,
\label{FE-1'}
\\
&&p_{12}(s)=\frac{\lambda_1}{a_{22}(f^{-1}(s))}p_{12}(f^{-1}(s))-\frac{a_{12}(f^{-1}(s))p_{11}(s)}{a_{22}(f^{-1}(s))},
\label{FE-2'}
\\
%&\tilde{p}_{22}(s)=\frac{\lambda_2}{a_{22}(f^{-1}(s))}\tilde{p}_{22}(f^{-1}(s))+\frac{\lambda_2}{a_{22}(f^{-1}(s))}-1,
&&\tilde{p}_{22}(s)=\frac{a_{22}(s)}{\lambda_2}\tilde{p}_{22}(f(s))+\frac{a_{22}(s)}{\lambda_2}-1,
\label{FE-3'}
\end{eqnarray}
where $\tilde{p}_{11}:=p_{11}-1$ and $\tilde{p}_{22}:=p_{22}-1$.
Note that equation (\ref{FE-2'}) is well defined since
$|a_{22}(s)|\ge |\lambda_2|-\eta>0$ by (\ref{FFbumpR2}) and
(\ref{FP-matrix}). In order to solve equation (\ref{FE-1'}), let
$E_1$ be a set of all $C^{0,\alpha}$ functions $\phi:I\to
\mathbb{R}$ such that $\phi(0)=0$ and
\begin{eqnarray*}
\|\phi\|_{E_1}:=\sup_{s\ne t\in I}
\frac{|\phi(s)-\phi(t)|}{|s-t|^\alpha}<\infty.
\end{eqnarray*}
Obviously, $E_1$ is a Banach space equipped with the norm $\|\cdot\|_{E_1}$ (cf. \cite{Evansbook98}).
Define an operator $T_1: E_1\to E_1$ by
\begin{eqnarray}
(T_1\phi)(s)=\frac{a_{11}(s)}{\lambda_1}\phi(f(s))+\frac{a_{11}(s)}{\lambda_1}-1,
\qquad \forall \phi\in E_1.
\label{T1}
\end{eqnarray}
One verifies that $T_1$ is indeed a self-mapping on $E_1$ since $a_{11}$
is $C^{0,\alpha}$ such that $a_{11}(0)=\lambda_1$. Furthermore, for all
$\phi,\psi\in E_1$, we have
%{\allowdisplaybreaks
\begin{eqnarray*}
&&\|(T_1\phi)-(T_1\psi)\|_{E_1}
\\
&&= |\lambda_1|^{-1}\frac{\big|a_{11}(s)(\phi-\psi)(f(s))-a_{11}(t)(\phi-\psi)(f(t))\big|}{|s-t|^\alpha}
\\
&&\le|\lambda_1|^{-1}\bigg(\frac{|a_{11}(s)||(\phi-\psi)(f(s))-(\phi-\psi)(f(t))|}{|s-t|^\alpha}
\\
    &&\hspace{2cm}+\frac{|a_{11}(s)-a_{11}(t)||(\phi-\psi)(f(t))|}{|s-t|^\alpha}\bigg)
\\
&&\le|\lambda_1|^{-1}\bigg(|a_{11}(s)|\frac{|(\phi-\psi)(f(s))-(\phi-\psi)(f(t))|}{|f(s)-f(t)|^\alpha}\frac{|f(s)-f(t)|^\alpha}{|s-t|^\alpha}
\\
    &&\hspace{2cm}+\frac{|a_{11}(s)-a_{11}(t)|}{|s-t|^\alpha}\frac{|(\phi-\psi)(f(t))|}{|f(t)|^\alpha}|f(t)|^\alpha\bigg)
\\
&&\le|\lambda_1|^{-1}\big\{(|\lambda_1|+\eta)^{1+\alpha}\|\phi-\psi\|_{E_1}+\delta \|a_{11}\|_{E_1}\|\phi-\psi\|_{E_1}\big\}
\\
&&\le (|\lambda_1|^\alpha+\varepsilon)\|\phi-\psi\|_{E_1},
\end{eqnarray*}
%\!\!\!}
where $\delta,\varepsilon>0$ are sufficiently small constants depending
on $I$. This implies that $T_1$ is a
contraction on the Banach space $E_1$ since $|\lambda_1|\in (0,1)$. Hence, by the Banach's contraction principle, $T_1$ has a unique fixed point in $E_1$, which is a solution of equation (\ref{FE-1'}).
Equation (\ref{FE-3'}) can be solved by the same method.

In order to solve equation (\ref{FE-2'}), let
$E_2$ be a space of all $C^{0,\alpha}$ functions $\phi:\mathbb{R}\to \mathbb{R}$ such that
$\phi(0)=0$, $\phi(s)=0$ for all $s\in \mathbb{R}\backslash J$ and
\begin{eqnarray*}
\|\phi\|_{E_2}:=\sup_{s\ne t\in J}
\frac{|\phi(s)-\phi(t)|}{|s-t|^\alpha}<\infty.
\end{eqnarray*}
Obviously, $E_2$ is a Banach space equipped with the norm $\|\cdot\|_{E_2}$.
Define an operator $T_2:E_2\to E_2$ by
\begin{eqnarray}
(T_2\phi)(s)=\frac{\lambda_1}{a_{22}(f^{-1}(s))}\phi(f^{-1}(s))-\frac{a_{12}(f^{-1}(s))p_{11}(s)}{a_{22}(f^{-1}(s))},
\qquad \forall \phi\in E_2.
\label{T2}
\end{eqnarray}
Here we know from (\ref{FFbumpR2}) and (\ref{FP-matrix}) that $a_{12}\in E_2$ and the solution $p_{11}$ of
equation (\ref{FE-1'}) obtained
above can be extended to the whole $\mathbb{R}$ with a compact support
as we did for $F$ before (\ref{FFbumpR2}). Therefore, $T_2$ is a self-mapping on $E_2$ since $f^{-1}$ is an
expansion by (\ref{fsbump}). Furthermore, for all $\phi,\psi\in E_2$,
\begin{eqnarray*}
&&\|(T_2\phi)-(T_2\psi)\|_{E_2}
\\
&&=\bigg|\frac{\lambda_1}{a_{22}(f^{-1}(s))}(\phi-\psi)(f^{-1}(s))
-\frac{\lambda_1}{a_{22}(f^{-1}(t))}(\phi-\psi)(f^{-1}(t))\bigg||s-t|^{-\alpha}
\\
&&\le\bigg|\frac{\lambda_1}{a_{22}(f^{-1}(s))}\bigg|\frac{|(\phi-\psi)(f^{-1}(s))-(\phi-\psi)(f^{-1}(t))|}{|f^{-1}(s)-f^{-1}(t)|^\alpha}
    \frac{|f^{-1}(s)-f^{-1}(t)|^\alpha}{|s-t|^{\alpha}}
\\
&&~~~~+ \bigg|\frac{\lambda_1}{a_{22}(f^{-1}(s))}
-\frac{\lambda_1}{a_{22}(f^{-1}(t))}\bigg|\,|s-t|^{-\alpha}\,\frac{|(\phi-\psi)(f^{-1}(t))|}{|f^{-1}(t)|^\alpha}|f^{-1}(t)|^\alpha
\\
&&\le (|\lambda_1|^{1-\alpha}|\lambda_2|^{-1}+\varepsilon)\|\phi-\psi\|_{E_2}.
\end{eqnarray*}
This implies that $T_2$ is a contraction on the Banach space $E_2$
since $|\lambda_1|^{1-\alpha}|\lambda_2|^{-1}\in (0,1)$ and $\varepsilon>0$ is small enough.
The last inequality of the above given formula holds because $|f^{-1}(t)|^\alpha\le \delta$ for all $t\in J$, where $\delta>0$
is sufficiently small provided $J$ is small enough.
Then the Banach's contraction principle gives the solution
of equation (\ref{FE-2'}) belonging to $E_2$. This completes the proof of Lemma \ref{lm-3FE}. \qquad$\Box$

%We will prove this lemma after we complete the proof of Lemma~\ref{????}
%%%
%\footnote{[ZWN130322] Fill up ???}
%%%
%.
By Lemma \ref{lm-3FE} we obtain the entries $p_{11}$, $p_{12}$, $p_{21}$ (which equals $0$)
and $p_{22}$ of the matrix $D\Psi(x_1,0)$, which satisfies equation (\ref{eqn-1}) for small $x_1$.
The same arguments carry over to the two equations given in
(\ref{eqn-2}). Namely, for (\ref{eqn-2}) it suffices to discuss the
equation
\begin{eqnarray}
D\Psi(F(0,x_2))\,DF(0,x_2)=\Lambda\, D\Psi(0,x_2)
\label{eqn-2}
\end{eqnarray}
and find a local $C^{1,\alpha}$ solution $\Psi$ of (\ref{eqn-2}) such that
$\Psi(O)=O$, $D\Psi(O)={\rm id}$ and $\pi_1\Psi(0,x_2)=0$. Similarly to
Lemma \ref{lm-3FE}, we can obtain the entries $q_{11}$, $q_{12}$ (which equals $0$), $q_{21}$ and $q_{22}$
of the matrix $D\Psi(0,x_2)$ which satisfies equation (\ref{eqn-2}) for small $x_2$.
Moreover, $q_{ij}$'s are $C^{0,\alpha}$ such that
\begin{eqnarray}
q_{11}(0)=q_{22}(0)=1
\quad {\rm and} \quad
q_{21}(0)=0.
\label{q1234}
\end{eqnarray}
Then we claim the following:
\begin{lm}
Let the $C^{0,\alpha}$ functions $p_{i,j}$'s and $q_{i,j}$'s are given above.
Then, there exists a $C^{1,\alpha}$ mapping $\Psi^*:\mathbb{R}^2\to\mathbb{R}^2$ such that
\begin{eqnarray}
\pi_1\Psi^*(0,x_2)=0, \quad \pi_2\Psi^*(x_1,0)=0
\label{flat-ceta}
\end{eqnarray}
and
\begin{eqnarray}
D\Psi^*(0,x_2)=
\left(\begin{array}{cc}
q_{11}(x_2) & 0
\\
q_{21}(x_2) & q_{22}(x_2)
\end{array}\right),
~~
D\Psi^*(x_1,0)=
\left(\begin{array}{cc}
p_{11}(x_1) & p_{12}(x_1)
\\
0 & p_{22}(x_1)
\end{array}\right)
\label{w-ext}
\end{eqnarray}
for all $x_1,x_2\in I$.
\label{lm-wa}
\end{lm}

{\bf Proof}. Our strategy of the proof is to use the Whitney extension theorem
(see Lemma \ref{lm-witn-ext} in the Appendix).
For the existence of $\pi_{1}\Psi^*$, in order to verify conditions (\ref{hmau}) and (\ref{whcond}) given just before Lemma \ref{lm-witn-ext},
we let $h_{ij}:(\{0\}\times I)\cup (I\times\{0\})\to \mathbb{R}$ ($i,j\in \{0,1\}$)
be functions defined by
\begin{eqnarray}
&h_{00}(x_1,x_2)=\left\{
\begin{array}{ll}
0, &~~ x_1=0,
\\
\int_{0}^{x_1}p_{11}(\tau)d\tau, &~~ x_2= 0,
\end{array}
\right.
\label{def-zeta1}
\end{eqnarray}
and
\begin{eqnarray}
h_{10}(x_1,x_2)=
\left\{
\begin{array}{ll}
q_{11}(x_2), &~~ x_1=0,
\\
p_{11}(x_1), &~~ x_2= 0,
\end{array}
\right.
\quad
h_{01}(x_1,x_2)=
\left\{
\begin{array}{ll}
0, &~~ x_1=0,
\\
p_{12}(x_1), &~~ x_2= 0.
\end{array}
\right.
\label{def-zeta2}
\end{eqnarray}
Obviously,
\begin{eqnarray}
|h_{ij}(x)|\le M, \qquad \forall x\in (\{0\}\times I)\cup (I\times\{0\}),
\label{wayyy1}
\end{eqnarray}
by the continuity of $p_{ij}$'s and $q_{ij}$'s.
Moreover, one checks that, for all $x,y\in (\{0\}\times I)\cup (I\times\{0\})$,
\begin{eqnarray}
&&h_{00}(x)=h_{00}(y)+h_{10}(y)(x_1-y_1)+h_{01}(y)(x_2-y_2)+R_{00}(x,y),
\label{eqn-zeta1}
\\
&&h_{10}(x)=h_{10}(y)+R_{10}(x,y),
\label{eqn-zeta2}
\\
&&h_{01}(x)=h_{01}(y)+R_{01}(x,y),
\label{eqn-zeta3}
\end{eqnarray}
where $R_{ij}:(\{0\}\times I)\cup (I\times\{0\})\to \mathbb{R}$ are functions such that
\begin{eqnarray}
|R_{ij}(x,y)|\le M\|x-y\|^{1+\alpha-(i+j)}, \quad \forall x,y\in (\{0\}\times I)\cup (I\times\{0\}),
\label{wayyy2}
\end{eqnarray}
for all $(i,j)=(0,0), (0,1)$ and $(1,0)$. In order to prove inequality (\ref{wayyy2}), note that $p_{ij}$'s
and $q_{ij}$'s are $C^{0,\alpha}$ functions satisfying (\ref{p1234}) and (\ref{q1234}) respectively. Then,
by (\ref{def-zeta1})-(\ref{def-zeta2}) and (\ref{eqn-zeta1})-(\ref{eqn-zeta3}), we have the following facts:

\noindent (i) For all $x,y\in \{0\}\times I$, i.e.,
$x_1=y_1=0$, we have
\begin{eqnarray*}
|R_{00}(x,y)|&=&|h_{00}(0,x_2)-h_{00}(0,y_2)-h_{01}(0,y_2)(x_2-y_2)|=0,
\\
|R_{10}(x,y)|&=&|h_{10}(0,x_2)-h_{10}(0,y_2)|=|q_{11}(x_2)-q_{11}(y_2)|\le L\|x-y\|^\alpha,
\\
|R_{01}(x,y)|&=&|h_{01}(0,x_2)-h_{01}(0,y_2)|=0.
\end{eqnarray*}

\noindent (ii) For all $x\in \{0\}\times I$ and all $y\in I\times \{0\}$, i.e.,
$x_1=y_2=0$, we have
\begin{eqnarray*}
|R_{00}(x,y)|
&=&|h_{00}(0,x_2)-h_{00}(y_1,0)+h_{10}(y_1,0)y_1-h_{01}(y_1,0)x_2|
~~~~~~~~~~
\\
&=&\Big|-\int_{0}^{y_1}p_{11}(\tau)d\tau+p_{11}(y_1)y_1-p_{12}(y_1)x_2\Big|
\\
&\le& |(p_{11}(\xi_1)-p_{11}(y_1))y_1+p_{12}(y_1)x_2|
\\
&\le& L(|y_1|^{1+\alpha}+|y_1^{\alpha}x_2|)\le 2L \max\{|y_1|^{1+\alpha},|x_2|^{1+\alpha}\}
\\
&=& 2L\|(0,x_2)-(y_1,0)\|^{1+\alpha}=2L\|x-y\|^{1+\alpha},
\end{eqnarray*}
where $\xi_1\in (0,y_1)$, and
\begin{eqnarray*}
|R_{10}(x,y)|
&=&|h_{10}(0,x_2)-h_{10}(y_1,0)|=|q_{11}(x_2)-p_{11}(y_1)|
~~~~~~~~~~~~~~~~~~~
\\
&\le&|q_{11}(x_2)-1|+|p_{11}(y_1)-1|\le|x_2|^\alpha+|y_1|^\alpha
\\
&\le& 2\|x-y\|^\alpha,
\\
|R_{01}(x,y)|
&=&|h_{01}(0,x_2)-h_{01}(y_1,0)|=|p_{12}(y_1)|\le |y_1|^\alpha
\\
&\le&\|x-y\|^\alpha.
\end{eqnarray*}

\noindent (iii) For all $x\in I\times \{0\}$ and all $y\in \{0\}\times I$, i.e.,
$x_2=y_1=0$, we have
\begin{eqnarray*}
|R_{00}(x,y)|
&=&|h_{00}(x_1,0)-h_{00}(0,y_2)-h_{10}(0,y_2)x_1+h_{01}(0,y_2)y_2|
~~~~~~~~~~
\\
&=&\Big|\int_{0}^{x_1}p_{11}(\tau)d\tau-q_{11}(y_2)x_1\Big|
\\
&=&|(p_{11}(\xi_2)-1)x_1-(q_{11}(y_2)-1)x_1|
\\
&\le& L(|x_1|^{1+\alpha}+|y_2^\alpha x_1|)\le 2L\max\{|x_1|^{1+\alpha}, |y_2|^{1+\alpha}\}
\\
&=& 2L\|x-y\|^{1+\alpha},
\end{eqnarray*}
where $\xi_2\in(0,x_1)$, and
\begin{eqnarray*}
|R_{10}(x,y)|
&=&|h_{10}(x_1,0)-h_{10}(0,y_2)|=|p_{11}(x_1)-q_{11}(y_2)|
~~~~~~~~~~~~~~~~~~~
\\
&\le&|x_1|^\alpha+|y_2|^\alpha\le 2\|x-y\|^\alpha,
\\
|R_{01}(x,y)|
&=&|h_{01}(x_1,0)-h_{01}(0,y_2)|=|p_{12}(x_1)|
\\
&\le&|x_1|^\alpha\le \|x-y\|^\alpha.
\end{eqnarray*}

\noindent (iv) For all $x,y\in I\times \{0\}$, i.e.,
$x_2=y_2=0$, we have
\begin{eqnarray*}
|R_{00}(x,y)|
&=&|h_{00}(x_1,0)-h_{00}(y_1,0)-h_{10}(y_1,0)(x_1-y_1)|
\\
&=&\Big|\int_{0}^{x_1}p_{11}(\tau)d\tau-\int_{0}^{y_1}p_{11}(\tau)d\tau-p_{11}(y_1)(x_1-y_1)\Big|
~~~~~~~~~~~
\\
&\le& |p_{11}(\xi_3)(x_1-y_1)-p_{11}(y_1)(x_1-y_1)|
\\
&\le& L|x_1-y_1|^{1+\alpha}=L\|x-y\|^{1+\alpha},
\end{eqnarray*}
where $\xi_3\in(x_1,y_1)$ or $(y_1,x_1)$, and
\begin{eqnarray*}
|R_{10}(x,y)|
&=&|h_{10}(x_1,0)-h_{10}(y_1,0)|=|p_{11}(x_1)-p_{11}(y_1)|
\\
&\le& L|x_1-y_1|^\alpha= L\|x-y\|^\alpha,
\\
|R_{01}(x,y)|
&=&|h_{01}(x_1,0)-h_{01}(y_1,0)|=|p_{12}(x_1)-p_{12}(y_1)|
~~~~~~~~~~~~~~~~~~~~
\\
&\le& L|x_1-y_1|^\alpha= L\|x-y\|^\alpha.
\end{eqnarray*}

According to (i)-(iv), we prove inequality (\ref{wayyy2}).
Hence, due to (\ref{wayyy1})-(\ref{wayyy2}), we see that (\ref{hmau}) and (\ref{whcond}) given in the appendix are satisfied.
It follows from Lemma \ref{lm-witn-ext} (see the appendix) that
there exists a $C^{1,\alpha}$ function $\pi_{1}\Psi^*:\mathbb{R}^2\to\mathbb{R}$ such that
$
\pi_1\Psi^*(0,x_2)=h_{00}(0,x_2)=0,
$
and
\begin{eqnarray*}
&&D(\pi_1\Psi^*)(0,x_2)=(h_{10}(0,x_2),h_{01}(0,x_2))=(q_{11}(x_2), 0),
\\
&&D(\pi_1\Psi^*)(x_1,0)=(h_{10}(x_1,0),h_{01}(x_1,0))=(p_{11}(x_1), p_{12}(x_1)).
\end{eqnarray*}
Similarly, we can find $\pi_{2}\Psi^*:\mathbb{R}^2\to \mathbb{R}$ such that
$
\pi_2\Psi^*(x_1,0)=0,
$
and
\begin{eqnarray*}
D(\pi_2\Psi^*)(0,x_2)=(q_{21}(x_2), q_{22}(x_2)),
~~~~~
D(\pi_2\Psi^*)(x_1,0)=(0, p_{22}(x_1)).
\end{eqnarray*}
This enables us to define $\Psi^*:\mathbb{R}^2\to \mathbb{R}^2$ by
$
\Psi^*:=(\pi_{1}\Psi^*,\pi_{2}\Psi^*),
$
which obviously satisfies (\ref{flat-ceta}) and (\ref{w-ext}).
This completes the proof of Lemma \ref{lm-wa}. \qquad$\Box$

Having Lemma 6, we see that $\Psi^*$ satisfies both
equations (\ref{eqn-1}) and (\ref{eqn-2}). Moreover,
$$
\Psi^*(O)=O, ~D\Psi^*(O)={\rm id}, ~\pi_1\Psi^*(0,x_2)=0 ~~{\rm and}~~ \pi_2\Psi^*(x_1,0)=0
$$
due to (\ref{p1234}) and (\ref{q1234})-(\ref{w-ext}).
Hence, by the above discussions for equations (\ref{eqn-1}) and (\ref{eqn-2}), we can prove
Lemma \ref{lm-flat}.

\begin{re}
{\rm
%Our Theorem~\ref{thm-C1-R2} only concerns the case of dimension 2.
If we consider the case of $\mathbb{R}^n$ ($n\ge 3$),
the dimension of the functional equations (\ref{FE-1})-(\ref{FE-3})
will be $\ge 2$ and therefore the Banach's contraction principle is not applicable because the operators $T_1$ and $T_2$,
defined in (\ref{T1}) and (\ref{T2}) respectively will be not a contraction any more in general.
This is the reason why we do not tackle the problem in $\mathbb{R}^n$ as Hartman and Belitskii did.
}
\end{re}

%%--------------------------------------
\section{Sharpness of the estimates on $\beta$}
\setcounter{equation}{0}

This section is devoted to the following result:

\begin{thm}
In the case of {\rm (\ref{sdcs})}, the estimates for $\beta$ given
in Theorem {\rm \ref{thm-C1-R2}} are sharp upper bounds.
\label{fa-sad1}
\end{thm}

{\bf Proof}:
We first consider the case that $|\lambda_1\lambda_2|\le 1$. Our strategy is to prove that
if $\Phi:U\to \mathbb{R}^2$ is a $C^{1,\beta}$ diffeomorphism which linearizes the
$C^{1,\alpha}$ mapping $F_*$ defined in (\ref{F-def}) then
$
\beta\le -\sigma(1-\sigma)^{-1}\alpha.
$

As indicated in the first paragraph of the proof of Theorem \ref{fa-con}, we suppose that $0<\lambda_1<1<\lambda_2$.
Fix a sufficiently small constant $\xi>0$ and put
\begin{eqnarray}
n_0(x_2):=(1-\sigma)^{-1}\log_{\lambda_1}(x_2/\xi)> 1
\label{n0}
\end{eqnarray}
for small $x_2>0$, where $\sigma:={\log\lambda_2}/{\log\lambda_1}<0$ is given such that $\lambda_2=\lambda_1^{\,\sigma}$.
We simply let $n_0$ denote $n_0(x_2)$ in what follows. Clearly, the
integer $[n_0]$, the greatest integer not exceeding $n_0$, is the largest integer in the set $\{n\in\mathbb{N}:\lambda_2^{n}x_2\le\lambda_1^{n}\xi\}$.
Observing (\ref{F-def}), we get
\begin{eqnarray}
\pi_1 F_*(x)\ge \lambda_1 x_1, ~~~~~ \forall  x\in \mathbb{R}^2,
\label{F1>r1}
\end{eqnarray}
because $\varrho(x)\ge 0$ and $u(x)\ge 0$.
Furthermore, one can check from (U1) and (U2) that $F_*$ satisfies (\ref{FwwbumpR2}) and (\ref{flat-Fww}) with $\widetilde{F}$ replacing with $F_*$
and consequently
\begin{eqnarray}
\lambda_1^{n}\xi\le\pi_1F_*^{n}(\xi,x_2)\le
M\lambda_1^{n}\xi, ~~~~~ \forall n\in \mathbb{N},
\label{lFl}
\end{eqnarray}
due to (\ref{F1F1'}) and (\ref{F1>r1}). One also computes directly that
\begin{eqnarray}
\pi_2 F_*^{n}(\xi,x_2)=\lambda_2^nx_2, \qquad \forall n\in\mathbb{N},
\label{pi2F*}
\end{eqnarray}
and that
%{\allowdisplaybreaks
\begin{align}
\!\!\!\pi_1 F_*^{n}(\xi,x_2)
&\ge \lambda_1^{n-[n_0]}\pi_1 F_*^{[n_0]}(\xi,x_2)
\nonumber\\
&=
\lambda_1^{n-[n_0]}\Big\{\lambda_1^{[n_0]}\xi+\sum_{i=1}^{[n_0]-1}\lambda_1^{[n_0]-i}\varrho(F_*^i(\xi,x_2))u(\Lambda F_*^i(\xi,x_2))
    (\lambda_2^ix_2)^{1+\alpha}
\nonumber\\
    &~~~~+\varrho(F_*^{[n_0]-1}(\xi,x_2))u(\Lambda F_*^{[n_0]-1}(\xi,x_2))(\lambda_2^{[n_0]}x_2)^{1+\alpha}\Big\}
\nonumber\\
&\ge \lambda_1^{n-[n_0]}\Big\{\lambda_1^{[n_0]}\xi\!+\!
     \varrho(F_*^{[n_0]-1}(\xi,x_2))u(\Lambda F_*^{[n_0]-1}(\xi,x_2))(\lambda_2^{[n_0]}x_2)^{1+\alpha}\Big\}
\label{pi1Fn0-1}
\end{align}
for small $x_2>0$ and for all $n\ge [n_0]$ by (\ref{F-def}) and (\ref{F1>r1}).

On the other hand, in view of (\ref{lFl}) and (\ref{pi2F*}) we have
\begin{eqnarray*}
\pi_1F_*^{[n_0]-1}(\xi,x_2)\le
M\lambda_1^{[n_0]-1}\xi,
\qquad
\pi_2F_*^{[n_0]-1}(\xi,x_2)=\lambda_2^{[n_0]-1}x_2\le \lambda_1^{[n_0]-1}\xi
\end{eqnarray*}
since $[n_0]-1$ is in the set $\{n\in\mathbb{N}:\lambda_2^{n}x_2\le\lambda_1^{n}\xi\}$. It implies that
$
\|F_*^{[n_0]-1}(\xi,x_2)\|\le M\lambda_1^{[n_0]-1}\xi\le M\xi
$
and therefore
\begin{eqnarray}
\varrho(F_*^{[n_0]-1}(\xi,x_2))=1
\label{bf=1}
\end{eqnarray}
since $\xi$ is small enough.
Furthermore, it follows from (U1) that
\begin{eqnarray}
u(\Lambda F_*^{[n_0]-1}(\xi,x_2))=u\big(\lambda_1 \pi_1F_*^{[n_0]-1}(\xi,x_2),\lambda_2^{[n_0]} x_2\big)=1
\label{u=1}
\end{eqnarray}
since $\lambda_2^{[n_0]}x_2\le\lambda_1^{[n_0]}\xi\le
\lambda_1\pi_1F_*^{[n_0]-1}(\xi,x_2)$ by (\ref{lFl}). Thus
(\ref{pi1Fn0-1})-(\ref{u=1}) yield
\begin{eqnarray}
\pi_1 F_*^{n}(\xi,x_2)
\ge\lambda_1^{n-[n_0]}\{\lambda_1^{[n_0]}\xi+(\lambda_2^{[n_0]}x_2)^{1+\alpha}\}
\ge \lambda_1^{n}\xi +M_1\lambda_1^nx_2^{(1-\sigma)^{-1}\alpha}.
\label{pi1Fn0-2}
\end{eqnarray}
by (\ref{n0}). Substituting $\lambda_2^{-n}\omega$ for $x_2$ in (\ref{pi1Fn0-2}) we obtain
\begin{eqnarray}
\pi_1 F_*^{n}(\xi,\lambda_2^{-n}\omega)-\lambda_1^{n}\xi \ge
M_2\lambda_1^{n(1-\sigma(1-\sigma)^{-1}\alpha)}, \qquad \forall n\ge
[n_0],
\label{tymmx}
\end{eqnarray}
where the constant $\omega>0$ is chosen to be small enough.

Now we are ready to prove the sharpness of the estimates of $\beta$ given in Theorem \ref{thm-C1-R2}.
Fix a constant $\beta\in (0,1]$ and suppose that $\Phi:U\to \mathbb{R}^2$ is a $C^{1,\beta}$ diffeomorphism such that
equation (\ref{schd-eqn}) with $F:=F_*$ holds near $O$.
According to (4.55) in \cite{WZWZ-JFA01}, without loss of generality, we may assume that
\begin{eqnarray}
\Phi(0,x_2)=(0,x_2),
~~~\Phi(x_1,0)=(x_1,0),
~~~{D\Phi}(x_1,0)={D\Phi}(0,x_2)={\rm id}.
\label{phi-phi'}
\end{eqnarray}
By (\ref{phi-phi'}), the Taylor expansion of $\pi_1\Phi$ at $(\xi,0)$ gives
$\pi_1\Phi(\xi,x_2)=\xi+O(x_2^{1+\beta})$. Substituting $\lambda_2^{-n}\omega$ for $x_2$, where $\omega>0$
is a small constant such that $(\xi, \lambda_2^{-n}\omega)\in U$ for all $n\in\mathbb{N}$,
we get
\begin{eqnarray}
\pi_1\Phi(\xi,\lambda_2^{-n}\omega)=\xi+O(\lambda_2^{-n(1+\beta)}), ~~~~~ {\rm as} ~ n\to \infty.
\label{phi-exp}
\end{eqnarray}
By (\ref{schd-eqn}) with $F:=F_*$ and (\ref{phi-exp}),
\begin{eqnarray}
\pi_1\Phi(F_*^{n}(\xi,\lambda_2^{-n}\omega))
&=&\lambda_1^{n}\pi_1\Phi(\xi,\lambda_2^{-n}\omega)=\lambda_1^{n}\xi+O((\lambda_1\lambda_2^{-(1+\beta)})^n)
\nonumber\\
&=&\lambda_1^{n}\xi+O(\lambda_1^{n(1-\sigma(1+\beta))}).
\label{Fphi-exp1}
\end{eqnarray}
On the other hand, the Taylor expansion of $\pi_1\Phi$ at $(0,\omega)$ gives
\begin{eqnarray}
\pi_1\Phi(F_*^{n}(\xi,\lambda_2^{-n}\omega))
&=&\pi_1\Phi(\pi_1F_*^{n}(\xi,\lambda_2^{-n}\omega),\omega)
\nonumber\\
&=&
\pi_1F_*^{n}(\xi,\lambda_2^{-n}\omega)+O(\{\pi_1F_*^{n}(\xi,\lambda_2^{-n}\omega)\}^{1+\beta})
\nonumber\\
&=&\pi_1F_*^{n}(\xi,\lambda_2^{-n}\omega)+O(\lambda_1^{n(1+\beta)})
\label{xnyyy}
\end{eqnarray}
by (\ref{lFl}).

Combining (\ref{Fphi-exp1}) with (\ref{xnyyy}), we get
\begin{eqnarray}
\pi_1F_*^{n}(\xi,\lambda_2^{-n}\omega)-\lambda_1^{n}\xi=O(\lambda_1^{n(1-\sigma(1+\beta))})+O(\lambda_1^{n(1+\beta)}).
\label{r3y}
\end{eqnarray}
Assume that $\beta:=\delta-\sigma(1-\sigma)^{-1}\alpha>0$ for an arbitrary constant $\delta>0$,
where $\sigma<0$ as indicated below (\ref{n0}).
Then one computes that
\begin{eqnarray*}
1-\sigma(1+\beta)=1-\sigma(1+\delta-\sigma(1-\sigma)^{-1}\alpha)
>1-\sigma(1-\sigma)^{-1}\alpha>0
\end{eqnarray*}
and
\begin{eqnarray*}
1+\beta=1+\delta-\sigma(1-\sigma)^{-1}\alpha>1-\sigma(1-\sigma)^{-1}\alpha>0.
\end{eqnarray*}
Therefore (\ref{r3y}) contradicts to (\ref{tymmx}) since $\lambda_1\in (0,1)$, which
means that
$$
\beta\le -\sigma(1-\sigma)^{-1}\alpha.
$$
The proof of the case $|\lambda_1\lambda_2|\le 1$ is completed.

Next we study the case that $|\lambda_1\lambda_2|> 1$. For the purpose, we give another diffeomorphism
\begin{eqnarray*}
G(x):=\left\{
\begin{array}{ll}
(\lambda_1^{-1}x_1, \lambda_2^{-1}x_2+\varrho(x){\tilde{u}}(\Lambda^{-1} x) p_\alpha(\lambda_1^{-1}x_1)), & x\in \mathbb{R}^2\backslash \{O\},
\\
\vspace{-0.3cm}
\\
O, & x=O,
\end{array}
\right.
%~~~~~
%\label{G-def}
\end{eqnarray*}
where $\tilde{u}(x_1,x_2):=u(x_2,x_1)$, and consider its inverse $G_*:=G^{-1}$.
Similarly to the case of $F_*$, one can verify that $G_*:\mathbb{R}^2\to \mathbb{R}^2$ is $C^{1,\alpha}$ such that
$DG_*(O)=\Lambda$. We conclude that if $\Phi:U\to \mathbb{R}^2$ is a $C^{1,\beta}$ diffeomorphism which linearizes $G_*$ then
$$
\beta\le (1-\sigma)^{-1}\alpha.
$$
Indeed, this can be proved by replacing $F_*$, $\lambda_1$, $\lambda_2$, $x_1$, $x_2$, $u$ and $\sigma$ with $G_*$,
$\lambda_2^{-1}$, $\lambda_1^{-1}$, $x_2$, $x_1$, $\tilde{u}$ and $\sigma^{-1}$ respectively in the above proof for $F_*$.
The proof is completed.
\qquad$\Box$

\begin{re}
{\rm Combining the results given in \cite{WZWZ-JFA01} (as mentioned
in the Introduction) with our Theorems
\ref{fa-con}-\ref{fa-sad1},
we can give sharp bounds for the exponents $\alpha$ and $\beta$ in
various cases. To sum up, we draw graphs of the sharp bounds of
$\beta$ as functions of $\alpha$ as follows:

%%%%%%%%%%%%%%%%%%%%%%%%%%%%%%%%%%%%%%%%%%%%%%%%%%%%%%%%%%%%%%%%%%%%%%%%%%%%%%
\vspace{0.2cm}
\begin{picture}(390,180)(0,0)

%***********************begin Figure 1
\!\!\!\!\!\!\!\!\!\!
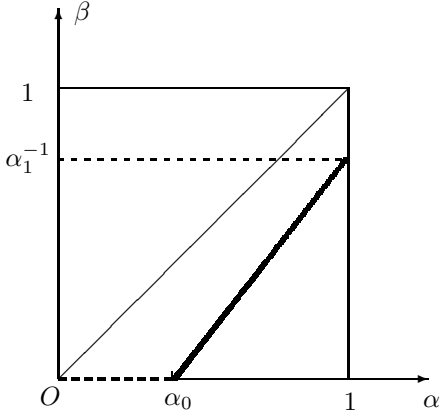
\put(10,0){\footnotesize Figure 1. In Poincar\'{e} domain: $\alpha_1\ge 1$}
%\put(45,-15){\footnotesize  when }
\put(71,37){\vector(1,0){96}}                                      %x_1-axis
\put(165,28){\makebox(2,1)[l]{\footnotesize $\alpha$}}
\put(27,37){\vector(0,1){140}}                                      %x_2-axis
\put(33,175){\makebox(2,1)[l]{\footnotesize $\beta$}}
\put(20,30){\makebox(2,1)[l]{{\footnotesize $O$}}}                  %origin

\put(27,147){\line(1,0){110}}                                       %horizontal

\put(136.7,37){\line(0,1){110}}                                     %vertical

\put(27,37){\line(1,1){109.7}}                                      %diagonal

%%%%%%%%%%%%%%%%%%%%%%%%%%%%%%%%%%%%%%%%%

\put(70,37){\line(0,1){3}}
\put(67,28){\makebox(2,1)[l]{\footnotesize $\alpha_0$}}

\put(135,28){\makebox(2,1)[l]{\footnotesize $1$}}

\put(13,145){\makebox(2,1)[l]{\footnotesize $1$}}

%\put(150,37){\line(0,1){3}}
%\put(147,28){\makebox(2,1)[l]{\footnotesize $\alpha_1$}}

%%%%%%%%%%%%%%%%%%%%%%%%%%%%%%%%%%%%%%%%%%

{\linethickness{0.5mm}\qbezier(71,37)(71,37)(135.7,120)}

%{\linethickness{0.5mm}\qbezier(27,37)(71,37)(71,37)}

\multiput(27,120)(5,0){22}{\line(1,0){2}}

{\linethickness{0.5mm}\multiput(27,37)(5,0){9}{\line(1,0){3}}}

\put(7,120){\makebox(2,1)[l]{\footnotesize $\alpha_1^{-1}$}}

%***********************end Figure 1

~~~~~~~~~~~~~~~~~~~~~~~~~~~~~~~~~~~~~~~~~~~~~~~~~~~ ~~~~~

%***********************begin Figure 2

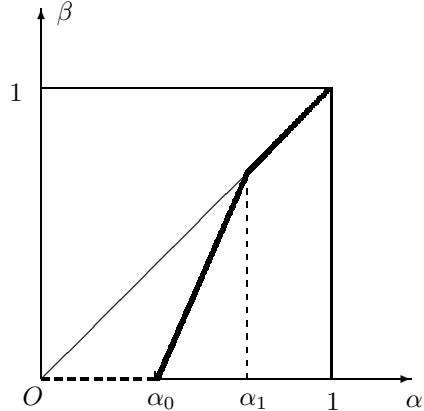
\put(5,0){\footnotesize Figure 2. In Poincar\'{e} domain: $0<\alpha_1<1$}
%\put(45,-15){\footnotesize  when }
\put(71,37){\vector(1,0){96}}                                      %x_1-axis
\put(165,28){\makebox(2,1)[l]{\footnotesize $\alpha$}}
\put(27,37){\vector(0,1){140}}                                      %x_2-axis
\put(33,175){\makebox(2,1)[l]{\footnotesize $\beta$}}
\put(20,30){\makebox(2,1)[l]{{\footnotesize $O$}}}                  %origin

\put(27,147){\line(1,0){110}}                                       %horizontal

\put(136.7,37){\line(0,1){110}}                                     %vertical

\put(27,37){\line(1,1){109.7}}                                        %diagonal

%%%%%%%%%%%%%%%%%%%%%%%%%%%%%%%%%%%%%%%%%

\put(70,37){\line(0,1){3}}
\put(67,28){\makebox(2,1)[l]{\footnotesize $\alpha_0$}}

\put(102,28){\makebox(2,1)[l]{\footnotesize $\alpha_1$}}

\put(135,28){\makebox(2,1)[l]{\footnotesize $1$}}

\put(15,145){\makebox(2,1)[l]{\footnotesize $1$}}

%%%%%%%%%%%%%%%%%%%%%%%%%%%%%%%%%%%%%%%%%%

{\linethickness{0.5mm}\qbezier(71,37)(71,37)(105,115)}

{\linethickness{0.5mm}\qbezier(105,115)(105,115)(135.7,146.5)}

%{\linethickness{0.5mm}\qbezier(27,37)(71,37)(71,37)}

{\linethickness{0.5mm}\multiput(27,37)(5,0){9}{\line(1,0){3}}}

\multiput(105,37)(0,5){15}{\line(0,1){2}}

%***********************end Figure 2

\end{picture}

\vspace{1cm}

\begin{picture}(390,180)(0,0)

\!\!\!\!\!\!\!\!\!\!
%***********************begin Figure 3

\put(10,0){\footnotesize Figure 3. In Siegel domain: $|\lambda_1\lambda_2|\le 1$}
%\put(45,-15){\footnotesize  }
\put(27,37){\vector(1,0){140}}                                      %x_1-axis
\put(165,28){\makebox(2,1)[l]{\footnotesize $\alpha$}}
\put(27,37){\vector(0,1){140}}                                      %x_2-axis
\put(33,175){\makebox(2,1)[l]{\footnotesize $\beta$}}
\put(20,30){\makebox(2,1)[l]{{\footnotesize $O$}}}                  %origin

\put(27,147){\line(1,0){110}}                                       %horizontal

\put(136.7,37){\line(0,1){110}}                                     %vertical

\put(27,37){\line(1,1){109.7}}                                      %diagonal

%%%%%%%%%%%%%%%%%%%%%%%%%%%%%%%%%%%%%%%%%

%\put(70,37){\line(0,1){3}}
%\put(67,28){\makebox(2,1)[l]{\footnotesize $\alpha_0$}}

\put(135,28){\makebox(2,1)[l]{\footnotesize $1$}}

\put(13,145){\makebox(2,1)[l]{\footnotesize $1$}}

%\put(150,37){\line(0,1){3}}
%\put(147,28){\makebox(2,1)[l]{\footnotesize $\alpha_1$}}

%%%%%%%%%%%%%%%%%%%%%%%%%%%%%%%%%%%%%%%%%%

{\linethickness{0.5mm}\qbezier(27,37)(27,37)(135.7,100)}

%{\linethickness{0.5mm}\qbezier(27,37)(71,37)(71,37)}

\multiput(27,100)(5,0){22}{\line(1,0){2}}

%{\linethickness{0.5mm}\multiput(27,37)(5,0){9}{\line(1,0){3}}}

\put(32,115){\makebox(2,1)[l]{\footnotesize
$-\sigma(1-\sigma)^{-1}$}}

\put(38,110){\vector(-1,-1){10}}

%***********************end Figure 3

~~~~~~~~~~~~~~~~~~~~~~~~~~~~~~~~~~~~~~~~~~~~~~~~~~~ ~~~~~

%***********************begin Figure 4

\put(5,0){\footnotesize Figure 4. In Siegel domain: : $|\lambda_1\lambda_2|>1$}
%\put(45,-15){\footnotesize  when }
\put(27,37){\vector(1,0){140}}                                      %x_1-axis
\put(165,28){\makebox(2,1)[l]{\footnotesize $\alpha$}}
\put(27,37){\vector(0,1){140}}                                      %x_2-axis
\put(33,175){\makebox(2,1)[l]{\footnotesize $\beta$}}
\put(20,30){\makebox(2,1)[l]{{\footnotesize $O$}}}                  %origin

\put(27,147){\line(1,0){110}}                                       %horizontal

\put(136.7,37){\line(0,1){110}}                                     %vertical

\put(27,37){\line(1,1){109.7}}                                      %diagonal

%%%%%%%%%%%%%%%%%%%%%%%%%%%%%%%%%%%%%%%%%

%\put(70,37){\line(0,1){3}}
%\put(67,28){\makebox(2,1)[l]{\footnotesize $\alpha_0$}}

\put(135,28){\makebox(2,1)[l]{\footnotesize $1$}}

\put(13,145){\makebox(2,1)[l]{\footnotesize $1$}}

%\put(150,37){\line(0,1){3}}
%\put(147,28){\makebox(2,1)[l]{\footnotesize $\alpha_1$}}

%%%%%%%%%%%%%%%%%%%%%%%%%%%%%%%%%%%%%%%%%%

{\linethickness{0.5mm}\qbezier(27,37)(27,37)(135.7,100)}

%{\linethickness{0.5mm}\qbezier(27,37)(71,37)(71,37)}

\multiput(27,100)(5,0){22}{\line(1,0){2}}

%{\linethickness{0.5mm}\multiput(27,37)(5,0){9}{\line(1,0){3}}}

\put(32,115){\makebox(2,1)[l]{\footnotesize $(1-\sigma)^{-1}$}}

\put(38,110){\vector(-1,-1){10}}

%***********************end Figure 4

\end{picture}

\vspace{0.8cm}

\noindent In Figures 1 and 2, the thick dashed lines mean that differentiable linearization at $O$ is permitted, but
there is counter example which cannot be $C^1$ linearized when $\alpha$ belongs to those areas.
}
\end{re}

\section*{Appendix: Whitney extension theorem}
Let ${\bf
j}:=(j_1,...,j_n)$, ${\bf k}:=(k_1,...,k_n)$, ${\bf
j}\,!:=j_1!\cdots j_n!$, $|\,{\bf j}\,|:=j_1+\cdots+j_n$ and $z^{\bf
j}:=z_1^{j_1}\cdots z_n^{j_n}$, where $n\in \mathbb{N}$, $j_i,k_i\ge
0$ are integers for all $i=1,...,n$ and $z:=(z_1,...,z_n)\in
{\mathbb{R}}^n$. For a constant $r\in \mathbb{N}$ and a function
$h:\Omega\to \mathbb{R}$, where $\Omega\subset \mathbb{R}^n$ is a
closed subset, we assume that there is a family of functions $h_{\bf
j}: \Omega\to \mathbb{R}$, $0< |\,{\bf j}\,|\le r$, such that
\begin{gather}
h_{\bf 0}=h, ~~~~~ h_{\bf j}(x)=\sum_{|{\bf j}+{\bf k}|\le
r}\frac{h_{{\bf j}+{\bf k}}(y)}{{\bf k}!}(x-y)^{\bf k}+R_{\bf
j}(x,y)
\label{hmau}
\tag{A.1}
\end{gather}
and
\begin{gather}
|h_{\bf j}(x)|\le M, ~~~~~ |R_{\bf j}(x,y)|\le
M\|x-y\|^{r+\alpha-|{\bf j}\,|}, ~~~~~ \forall x,y\in \Omega,
\tag{A.2}
\label{whcond}
\end{gather}
with $\alpha\in (0,1]$. Then, in view of Theorem 4 given in
\cite[p.177]{stein-book1970} (see also \cite{Whit-TAMS34}), we have
the following lemma:
\begin{lm}
Let the functions $h:\Omega\to \mathbb{R}$ and $h_{\bf j}: \Omega\to
\mathbb{R}$, $0\le |\,{\bf j}\,|\le r$, be given above such that
{\rm (\ref{hmau})} and {\rm (\ref{whcond})} hold. Then there exists
a $C^{r,\alpha}$ function $H:\mathbb{R}^n\to \mathbb{R}$ such that
\begin{eqnarray*}
H|_{\Omega}=h, ~~~~~~~ D^{\bf j}H|_{\Omega}=h_{\bf j}, ~~~~~ 0<
|\,{\bf j}\,|\le r,
\end{eqnarray*}
where $D^{\bf j}H:=\partial^{|{\bf j}\,|} H(x_1,...,x_n)/\partial
x_{j_1}\cdots\partial x_{j_n}$.
\label{lm-witn-ext}
\end{lm}

%%%%%%%%%%%%%%%%%%%%%%%%%%%%%%%%%%%%%%%%%%%%%
{\footnotesize


\begin{thebibliography}{79}



\bibitem{AraBeliZhu-book96} S. Kh. Aranson, G. R. Belitskii and E. V. Zhuzhoma, {\it Introduction to the Qualitative
Theory of Dynamical Systems on Surfaces}, Amer. Math. Soc., Providence, 1996.

\bibitem{Arnold83} V. I. Arnold, {\it Geometric Methods in the Theory of Ordinary Differential Equations},
                   Springer-Verlag, New York, 1983.

\bibitem{AAIS-book94} V. I. Arnold, V. S. Afrajmovich, Yu. S. Il'yashenko and L. P. Shil'nikov, {\it Dynamical Systems
V: Bifurcation Theory and Catastrophe Theory}, Encyclop.
%%Encyclopaedia
Math. Sci. {\bf 5}, Springer, Berlin, 1994.

\bibitem{Beli-RMS78} G. R. Belitskii, Equivalence and normal forms of germs of smooth mappings,
                     {\it Russian Math. Surveys} {\bf 33} (1978), 107-177.


\bibitem{Bronkopa-book94} I. U. Bronstein and A. Ya. Kopanski,
                          {\it Smooth Invariant Manifolds and Normal Forms}, World Scientific, River Edge, NJ, 1994.

\bibitem{Chap-TMIS02} M. Chaperon, Invariant manifolds revisited, {\it Tr. Mat. Inst. Steklova}
{\bf 236} (2002), 428-446,
Dedicated to the 80th anniversary of Academician Evgenii Frolovich Mishchenko, (Russian) (Suzdal, 2000).

\bibitem{ChHaTan-JDE97} X.-Y. Chen, J. K. Hale and B. Tan, Invariant foliations for $C^1$ semigroups in Banach spaces,
                        {\it J. Differential Equations} {\bf 139} (1997), 283-318.

\bibitem{DBo-JDE89} B. Deng, The Sil'nikov problem, exponential expansion, strong $\lambda$-lemma, $C^1$-linearization, and homoclinic bifurcation,
                      {\it J. Differential Equations} {\bf 79} (1989), 189-231.

\bibitem{Evansbook98} L. C. Evans, {\it Partial Differential Equations}, Amer. Math. Soc., Providence, 1998.


\bibitem{El-JFA01} M. S. ElBialy, Local contractions of Banach spaces and
spectral gap conditions, {\it J. Funct. Anal.} {\bf 182} (2001),
108-150.


\bibitem{FMT-Ann07} M. Field, I. Melbourne and A. T\"{o}r\"{o}k, Stability of mixing and rapid mixing for hyperbolic flows.
                    {\it Ann. Math.} {\bf 166} (2007), 269-291.

\bibitem{GuyHassRay-DCDS03} M. Guysinsky, B. Hasselblatt and V. Rayskin, Differentiability of the Hartman-Grobman linearization,
                             {\it Discr. Contin. Dyn. Syst.} {\bf 9} (2003), 979-984

\bibitem{Hart60} P. Hartman, On local homeomorphisms of Euclidean spaces, {\it Bol. Soc. Mat. Mexicana} {\bf 5} (1960), 220-241.

\bibitem{Hart64} P. Hartman, {\it Ordinary Differential Equations}, John Wiley \& Sons, New York, 1964.

\bibitem{HPS-book77} M. Hirsch, M. Shub and C. Pugh, {\it Invariant Manifolds}, Lecture Notes in Math. {\bf 583}, Springer, 1977.

\bibitem{HollMelb-07} M. Holland and I. Melbourne, Central limit theorems and invariance principles for Lorenz attractors,
                      {\it J. London Math. Soc.} {\bf 76} (2007), 345-364.

\bibitem{HomHow-CMP99} A. J. Homburg and H. Weiss, A geometric criterion for positive topological entropy
                       II: homoclinic tangencies,  {\it Commun. Math. Phys.} {\bf 208} (1999), 267-273.

\bibitem{Kuczma1968} M. Kuczma, {\it Functional Equations in a Single Variable},  Polish Scientific Publ., Warsaw, 1968.


%\bibitem{GGS-book10} F. Gazzola, H.-C. Grunau and G. Sweers, {\it Polyharmonic boundary value problems},
%                     Lecture Notes in Math. {\bf 1991}, Springer-Verlag, Berlin, 2010.

\bibitem{Poincare} H. Poincar\'{e}, Sur le probl\'{e}me des trois corps et les \'{e}quations de la dyanamique,
                     {\it Acta. Math.} {\bf 13} (1890), 1-270.

\bibitem{pugh-AJM69} C. Pugh, On a theorem of P. Hartman, {\it Amer. J. Math.} {\bf 91} (1969), 363-367.

\bibitem{Ray-JDE98} V. Rayskin, $\alpha$-H\"{o}lder linearization, {\it J. Differential Equ.}
                   {\bf 147} (1998), 271-284.

\bibitem{R-S-JDE04} H. M. Rodrigues and J. Sol${\rm \grave{a}}$-Morales, Linearization of class $C^1$ for contractions on
Banach spaces, {\it J. Differential Equations} {\bf 201} (2004),
351-382.

\bibitem{R-S-JDDE04} H. M. Rodrigues and J. Sol${\rm \grave{a}}$-Morales, Smooth linearization for a saddle on Banach spaces,
{\it J. Dyn. Differential Equations} {\bf 16} (2004),
767-793.

\bibitem{Sell-AJM85} G. R. Sell, Smooth linearization near a fixed point, {\it Amer. J. Math.} {\bf 107} (1985), 1035-1091.


\bibitem{stein-book1970} E. M. Stein, {\it Singular Integrals and Differentiability Properties of Functions,}
                         Princeton University Press, Princeton, NJ, 1970.

\bibitem{ster-Duke} S. Sternberg, Local $C^n$ transformations of the real line, {\it Duke Math. J.} {\bf24}(1957), 97-102.

\bibitem{ster57} S. Sternberg, Local contractions and a theorem of Poincar\'{e},
                {\it Amer. J. Math.} {\bf 79} (1957), 809-824.


\bibitem{ster58} S. Sternberg, On the structure of local homeomorphisms of Euclidean $n$-space,
                {\it Amer. J. Math.} {\bf 80} (1958), 623-631.

\bibitem{Stowe-JDE86} D. Stowe, Linearization in two dimensions, {\it J. Differential Equations} {\bf 63} (1986), 183-226.

\bibitem{Stri-JDE90} S. van Strien, Smooth linearization of hyperbolic fixed points without resonance conditions,
                    {\it J. Differential Equations} {\bf 85} (1990), 66-90.

\bibitem{Tan-JDE00} B. Tan, ${\sigma}$-H\"{o}lder continuous linearization near hyperbolic
                     fixed points in $\mathbb{R}^n$, {\it J. Differential Equations} {\bf 162} (2000), 251-269.

\bibitem{Whit-TAMS34} H. Whitney, Analytic extensions of differentiable functions defined in closed sets, {\it Trans. Amer. Math. Soc.}
                      {\bf 36} (1934), 63-89.

\bibitem{ZZJZ-JMAA12} W. M. Zhang, Y. Y. Zeng, W. Jarczyk and W. N. Zhang,
           Local $C^1$ stability versus global $C^1$ unstability for iterative roots,
           {\it J. Math. Anal. Appl.} {\bf 386} (2012), 75-82.

\bibitem{WZWZ-JFA01} W. M. Zhang and W. N. Zhang, $C^1$ linearization for planar contractions,
                   {\it J. Funct. Anal.} {\bf 260} (2011), 2043-2063.




\end{thebibliography}
\end{document}